\documentclass[a4paper,10pt]{article}
\usepackage{amssymb}
\usepackage{amsmath}
\usepackage{graphicx}
\usepackage{color}
\usepackage[usenames,dvipsnames]{xcolor}
\usepackage{soul}

\usepackage{pst-plot}
\usepackage{tikz}
\usetikzlibrary{arrows}

\DeclareGraphicsExtensions{.pdf,.png,.jpg}

\newcommand{\bigo}{{\cal O}}
\newcommand{\eps}{\varepsilon}
\newcommand{\dom}[1]{\mbox{\rm dom}\,#1}
\newcommand{\doa}[1]{\mbox{\rm doa}\,#1}

\def\complex{\mathbb{C}}
\def\reals{\mathbb{R}}
\def\nat{\mathbb{N}}
\def\integers{\mathbb{Z}}
\def\W{{\cal W}}
\def\D{\mathcal{D}}

\newcommand{\im}[1]{{\rm Im}\left(#1\right)}
\newcommand{\re}[1]{{\rm Re}\left(#1\right)}

\def\qed{{\hfill{\vrule height5pt width3pt depth0pt}\medskip}}

\newcommand{\dist}{{\rm dist }}
\newcommand{\PPC}[1]{{\rm PPC}$(#1)$}
\newcommand{\PPCS}[1]{{\rm PPC}$(#1)$}
\newcommand{\NPC}[1]{{\rm NPC}$(#1)$}
\newcommand{\ro}{{\varrho }}

\newtheorem{theorem}{Theorem}
\newtheorem{lemma}[theorem]{Lemma}

\newtheorem{definition}{Definition}
\newtheorem{remark}[theorem]{Remark}

\newtheorem{ex}{Example}

\renewcommand{\phi}{\varphi}

\title{On the Petras algorithm for verified integration of piecewise analytic functions}
\author{Ma{\l}gorzata Moczurad, Piotr Zgliczy\'nski}

\begin{document}

\maketitle

\begin{abstract}

We consider the algorithm for verified integration of piecewise analytic functions presented in Petras' paper~\cite{P02}.
The analysis of the algorithm contained in that paper is limited to a narrow class of functions and gives upper bounds only. We present an estimation of the complexity (measured by a number of evaluations of an integrand) of the algorithm, both upper and lower bounds, for a wider class of functions. We show
examples with complexity $\Theta(|\ln\eps|/\eps^{p-1})$, for any $p >1$, where $\eps$ is the desired accuracy of the computed integral.
\end{abstract}


\section{Introduction}
In the paper we discuss the complexity of the Petras' algorithm for verified integration of piecewise analytic functions \cite{P02}. Reader should be warned that  we measure the complexity by a number of evaluations of an integrand. Therefore
 we omit the complexity of subroutines in the algorithm, which the full complexity analysis of the algorithm should take into account. We should mention the proper treatment of these issues will require a definition of computable analytic functions. Such definitions
 exist in the literature (see for example~\cite{GaertnerHotz12, Hoeven05} and the references given there), however in the present work we skip these issues and concentrate on the geometric aspects of the problem.

Our task is, given $f$ a piecewise analytic function on $[a,b]$ and $\eps>0$, to find a number $\mathcal{I}$, such that
\begin{equation*}
  \left|\int_a^b f(x) dx - \mathcal{I} \right| \leqslant \eps.  
\end{equation*}
If $f$ is analytic, then there are algorithms for integration
(based for example on Gauss-Legendre quadrature)
with the apparent complexity $\bigo(|\ln \eps|)$, where $\eps$ is the desired accuracy. It should be stressed that the constant in $\bigo(|\ln \eps|)$ depends
on the function $f$, mainly on the shape of its domain of analyticity; see~\cite{P98} for a discussion of optimal quadratures depending on the domain of analyticity.

 However, a disturbance of analyticity at some points makes the convergence rate deteriorate.
For the sake of discussion, let us call them breakpoints or singularities.

Knowing the location of singularities is not enough; we cannot simply partition the interval at those points and apply a standard algorithm for analytic functions to smaller intervals. The problem is that algorithms require
the function to be analytic on the interval of integration {\em and} its neighbourhood (disk or ellipse containing the interval in the complex plane).

Our case is (seemingly) even more difficult: we know neither the number, nor the location of the breakpoints; we might not even know whether they exist.

We consider Petras' algorithm for  verified integration of piecewise analytic functions.
The analysis of the algorithm contained in~\cite{P02} is not sufficient, because the class of functions considered
is unnaturally restricted and only an upper bound for this class is given.

The main result of the paper is a more sophisticated estimation of the complexity
of the algorithm for a wider class of functions. Our results still do not cover the whole range of (piecewise) analytic functions,
but only those that satisfy Petras-type conditions (PPC and NPC) of the order $p$ (see Definitions~\ref{def:petras-routines} and~\ref{def:petras-negative}).
The class of functions considered by Petras in \cite{P02} satisfies PPC condition of order $p \leqslant 1$.

The main idea explored in our paper is that the complexity  depends mainly on the region of analyticity. The difference between a ``simple'' and a ``difficult'' function is not that the former is analytic, while the latter is not. Even an analytic function might be hard for integration, if it has singularities very close to the real axis.  On the other hand, a piecewise analytic function might be relatively simple, if the region of analyticity around breakpoints is wide ($p \leqslant 1$ in the Definition~\ref{def:region-D-p}).

Originally in~\cite{P02} Petras showed that the complexity
is $\bigo\left(\ln^2 \eps \right)$ for a class of functions satisfying a condition, which in our paper is PPC condition of order $p=1$. We show that the complexity depends on the order $p$ of PPC and NPC conditions that an integrated function satisfies, i.e.\ if $p > 1$, then the complexity is ${\Theta}\left(|\ln \eps|/\eps^{p-1} \right)$, while for $p\leqslant 1$ the complexity is indeed ${\Theta}\left(\ln^2 \eps \right)$. Moreover, we show examples of functions (see Section~\ref{sec:examples}) for which the complexity scales as $|\ln \eps|/\eps^{p-1}$ for any $p>1$.

Notice that the complexity analysis focuses on the phenomena occurring in the close neighbourhood of singular points, since the most significant increase in the number of evaluations of the integrand takes place there. The number of evaluations
generated in the last step of the algorithm (i.e.\ near breakpoints) is comparable to the number of evaluations generated in all previous steps. Moreover, properties considered do not occur globally, but only in the proximity of the breakpoints.

 The paper is organized in the following way. Section~\ref{sec:notations}  contains basic definitions and notations. In Section~\ref{sec:Petras-algorithm} we present the Petras' algorithm. In Section~\ref{sec:alg-analysis} we introduce main tools for the analysis of the Petras' algorithm; in particular the Petras-type conditions PPC and NPC are given there. In Sections~\ref{sec:lower-bnd-estm} and~\ref{sec:upper-bound} we show lower and upper bounds for the complexity of the Petras' algorithm. Section~\ref{sec:complexity} contains the main theorem of the paper.

\section{Notations and core definitions}\label{sec:notations}

As usual, by $\nat$, $\integers$, $\reals$, $\complex$ we denote the sets of natural (including 0), integer, real and complex numbers, respectively. By $\nat_+$ and $\reals_+$ we denote
the set of positive natural and real numbers. We use $\overline{A}$ to denote the closure of a set $A$.
We will use $|z|$ to denote the absolute value of $z \in \complex$ or $z \in \reals$.
For a point $x \in \complex$  and a set $Z \subset \complex$ we define
\begin{equation*}
    \mbox{dist}(x,Z) =  \inf_{y \in Z} |x-y|. 
\end{equation*}

\begin{definition}
Assume that  $a,b \in \reals$ and $a<b$.
\begin{enumerate}
\item We say that a function $f: \complex\supset \dom{f}\to \complex$ is a {\em piecewise analytic function on $[a,b]$} if there exist open, pairwise disjoint sets $\D_j\subset\complex$ such that
\begin{enumerate}
\item $f$ is analytic on $\D_j$ and
$$\dom{f}\cap \{x+i y: a\leqslant x \leqslant b\} \subset \overline{\bigcup_{j =1}^{\infty} \D_j},$$

\item $\D_j\cap [a,b] = (a_j, b_j)$ and
$$[a, b] = \overline{\bigcup_{j = 1}^{\infty} (a_j, b_j)}$$
\end{enumerate}
\item The {\em domain of analyticity of $f$} is
$$\doa{f} = \bigcup_{j =1}^{\infty} \D_j.$$
\end{enumerate}
\end{definition}

\begin{ex}
Some examples of piecewise analytic functions on $[-1,1]$:
\begin{itemize}
\item all analytic functions, such that their domain contains $[-1,1]$,
\item the function $f(z) = \exp(-1/z^2)$,
\item the function $f(z) = |\sin(1/z)|$.
\end{itemize}
\end{ex}

Let $S = \{s_1, \ldots, s_m\}\subset \reals$, where
$a \leqslant s_1 < \ldots < s_m \leqslant b$, be a set of points such that the function $f$ is analytic in each open interval $(s_i, s_{i+1})$ for $i = 1, \ldots, m-1$.

\begin{figure}[htb]
\centering

\begin{tikzpicture}

\draw[draw=black ] (1,0) parabola (2.5,1.15);
\filldraw[fill=gray!20!white, draw=gray!20!white ] (1,0) parabola (2.5,1.15) -- (2.5,0);
\filldraw[fill=gray!20!white, draw=gray!20!white ] (1,0) parabola (2.5,-1.15) -- (2.5,0);
\draw (1 cm,1pt) -- (1 cm,-1pt) node[anchor=north] {\small $s_{i-1}$};
\filldraw[fill=gray!20!white, draw=gray!20!white ] (4,0) parabola (2.5,1.15) -- (2.5,0);
\filldraw[fill=gray!20!white, draw=gray!20!white ] (4,0) parabola (2.5,-1.15) -- (2.5,0);
\draw (4 cm,1pt) -- (4 cm,-1pt) node[anchor=north] {\small $s_i$};
\draw[draw=black ] (1,0) parabola (2.5,1.15);
\draw[draw=black ] (1,0) parabola (2.5,-1.15);
\draw[draw=black ] (4,0) parabola (2.5,1.15);
\draw[draw=black ] (4,0) parabola (2.5,-1.15);

\filldraw[fill=gray!20!white, draw=gray!20!white ] (4,0) parabola (5.0,0.6) -- (5.0,0);
\filldraw[fill=gray!20!white, draw=gray!20!white ] (4,0) parabola (5.0,-0.6) -- (5.0,0);
\draw (6 cm,1pt) -- (6 cm,-1pt) node[anchor=north] {\small $s_{i+1}$};
\filldraw[fill=gray!20!white, draw=gray!20!white ] (6,0) parabola (5.0,0.6) -- (5.0,0);
\filldraw[fill=gray!20!white, draw=gray!20!white ] (6,0) parabola (5.0,-0.6) -- (5.0,0);
\draw[draw=black ] (4,0) parabola (5.0, 0.6);
\draw[draw=black ] (4,0) parabola (5.0, -0.6);
\draw[draw=black ] (6,0) parabola (5.0, 0.6);
\draw[draw=black ] (6,0) parabola (5.0, -0.6);

\draw[color = black] (0,0) -- (7,0);
\end{tikzpicture}
\caption{Region $\D_{\gamma, S}^p$ between $s_{i-1}$ and $s_{i+1}$ with $p > 1$.}\label{fig:D-p>1}
\end{figure}
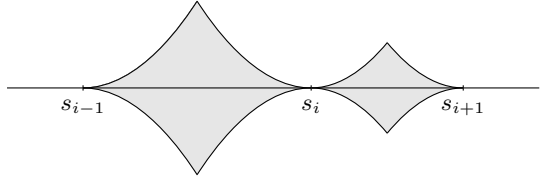

\begin{definition}\label{def:region-D-p}
For a given $\gamma$, $p > 0$ and a set $S$ define a region (see Figure~\ref{fig:D-p>1})
\begin{equation*} 
\D^p_{\gamma, S} := \{x + iy: |y|\leqslant \gamma\cdot \dist(x, S)^p\}.
\end{equation*}
\end{definition}

When $S$ is a singleton and its element is clear from the context, we omit the subscript $S$ and write $\D^p_{\gamma}$.
Notice that $\mathcal{V}_{\gamma, S}$ used in~\cite{P02} is a special case of $\D^p_{\gamma, S}$, i.e.\ $\mathcal{V}_{\gamma, S} = \D^1_{\gamma, S}$.

As usual we use
${\cal Q}_{n,[\alpha, \beta]}[f]$ to indicate an $n$-point quadrature formula for the evaluation of an integral of $f$ over $[\alpha,\beta]$.
$R_{n,[\alpha, \beta]}[f]$ is the error of the quadrature, i.e.
$$R_{n,[\alpha, \beta]}[f] = \left| \int_{\alpha}^{\beta} f(x)dx - {\cal Q}_{n,[\alpha, \beta]}[f]\right|. $$

\begin{figure}[htb]
\centering
\begin{tikzpicture}
\draw[dashed, fill = gray!20!white, draw =black] (1.1,-0.6) rectangle (3.4,0.6);
\draw[->] (0,0) -- (4,0);
\draw[->] (0.3,-1) -- (0.3,1.5);

\draw (3.0,1pt) -- (3.0,-1pt) node[anchor=south] {\small $\beta$};
\draw (1.5,1pt) -- (1.5,-1pt) node[anchor=south] {\small \!\!\! $\alpha$};
\draw (2.25,1pt) -- (2.25,-1pt) node[anchor=north] {\small $\frac{\alpha+\beta}{2}$};
\draw (0.2,-23pt) -- (0.4,-23pt);
\draw (0.2, -23pt) -- (0.25,-23pt) node[anchor=east] {\small \mbox{} $-\frac{d}{2}$};
\draw (0.2,23pt) -- (0.4,23pt);
\draw (0.2, 23pt) -- (0.25,23pt) node[anchor=east] {\small \mbox{} $\frac{d}{2}$};

\draw (3.0,0pt) -- (3.0,-30pt);
\draw (1.5,0pt) -- (1.5,-30pt);
\draw[->] (2.0, -30pt) -- (3.0, -30pt) node[anchor=north east] {\small $d = \beta - \alpha$};
\draw[->] (2.0, -30pt) -- (1.5, -30pt);

\end{tikzpicture}
\caption{Rectangle $\ro(\alpha, \beta,A,B) = \{x + iy: \frac{\alpha + \beta}{2} - A\frac{d}{2} \leqslant x \leqslant \frac{\alpha + \beta}{2} + A\frac{d}{2},\ -B\frac{d}{2}\leqslant y \leqslant B\frac{d}{2} \}$ }\label{fig:rectangle}
\end{figure}
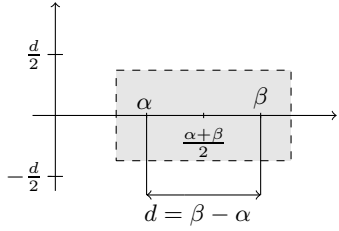

\begin{definition}\label{def:rho-rectangle}
For fixed $A> 1$, $B>0$, let the rectangle $\ro(\alpha, \beta, A, B)$ in the complex plane (see Figure~\ref{fig:rectangle}) be defined as
\begin{equation*}
\ro(\alpha, \beta,A,B) = \left\{ x + iy: \left|x - \frac{\beta + \alpha}{2}\right| \leqslant A \frac{\beta - \alpha}{2}, |y|\leqslant B\frac{\beta - \alpha}{2}\right\}.
\end{equation*}
When parameters $A$, $B$ are clear from the context, we omit them and write $\ro(\alpha, \beta)$.
\end{definition}

\begin{definition}
For fixed $A>1$, $B>0$, let $m(f; \alpha, \beta, A, B)$ be given by
$$m(f;\alpha, \beta, A, B) = \sup \{|f(z)|:\ z\in \ro(\alpha, \beta,A,B)\ \cap\ \dom{f}\}.$$
When parameters $A$, $B$ are clear from the context, we omit them and write $m(f;\alpha, \beta)$.
\end{definition}

Originally in~\cite{P02} the functional $m$ has a different meaning: it is a functional (or perhaps an algorithm) which returns a value greater or equal to $$\sup \{|f(z)|:\ z\in \ro(\alpha, \beta,A,B)\ \cap\ \dom{f}\},$$ which in this paper is realized by a function {\tt ComplexBound} (see the beginning of Section~\ref{sec:Petras-algorithm}).

For the estimation of error of the quadrature we need to see what is the size of the largest ellipse with the foci at $(\alpha,0)$ and $(\beta,0)$
which is contained in $\ro(\alpha,\beta, A, B)$.

\begin{lemma} \label{lem:ellipseAB}
Let $A>1$. If an ellipse with foci at $F_1 = (\alpha,0)$, $F_2 = (\beta, 0)$, major $\mathfrak{a} = A(\beta - \alpha)/2$ and minor semi-axis $\mathfrak{b} = B(\beta - \alpha)/2$ is inscribed (i.e.\ tangent to the edges) in $\ro(\alpha, \beta, A, B)$, then $B = \sqrt{A^2 - 1}$.
\end{lemma}
\begin{proof}
By definition of an ellipse
the distance from the center of the ellipse to the focal point is $\sqrt{\mathfrak{a}^2 - \mathfrak{b}^2}$, thus
\begin{equation*}
\frac{\beta - \alpha}{2} = \sqrt{\left(\frac{A(\beta - \alpha)}{2} \right)^2 - \left(\frac{B(\beta - \alpha)}{2} \right)^2}.
\end{equation*}
This gives
$$B = \sqrt{A^2 - 1}.$$
\qed
\end{proof}

By the above lemma parameters $A$ and $B$ of $\ro$ are not independent if we want to inscribe an ellipse into $\ro(\alpha,\beta, A, B)$. Therefore, in the sequel, we will use $B = \sqrt{A^2 - 1}$, only.

For $[\alpha, \beta] = [-1,1]$ we have $\mathfrak{a} = A$ and $\mathfrak{b} = B$. In particular for $A=\frac{5}{4}$ we have $B=\frac{3}{4}$ and the largest ellipse contained in $\ro(-1,1,A,B)$ has $\mathfrak{c}=\mathfrak{a}+\mathfrak{b}=2$.


\begin{remark}
\label{rem:gauss-cc-errors}
Let $A = \frac{5}{4}$ and $B = \frac{3}{4}$.  Then
\begin{eqnarray}
R_{n,[\alpha, \beta]}[f] & \leqslant & 2\cdot 4^{-n}\cdot (\beta - \alpha)m(f; \alpha, \beta) \quad\mbox{ (Gaussian quadrature)}
\label{eq:gauss-error}
\\[1ex]
R_{n,[\alpha, \beta]}[f] & \leqslant & 3\cdot 2^{-n}\cdot (\beta - \alpha)m(f; \alpha, \beta) \quad\mbox{ (Clenshaw-Curtis quadrature)}\label{eq:c-c-error}
\end{eqnarray}

for functions analytic on $\ro(\alpha,\beta, A, B)$.
\end{remark}
\begin{proof}
This not a real proof, but an explanation of the denominators in these formulas. For references regarding the constants used see~\cite{Br} and~\cite{P98-2}.

The map $$T(z)= \frac{\alpha+ \beta}{2} + \frac{\beta-\alpha}{2}z$$ is an affine isomorphism, such that $T(\ro(-1,1,A,B))=\ro(\alpha,\beta,A,B)$. The formulas
for Gauss or Clenshaw-Curtis quadratures and their errors are transported from $[-1,1]$ to $[\alpha,\beta]$ by this map.

On the normalized interval $[-1,1]$ the error of Gauss quadrature contains the term $1/ \mathfrak{c}^{2n}$, where $\mathfrak{c}$ is the sum of major and minor semi-axes.
 The transformation of $[\alpha,\beta]$ onto $[-1,1]$
multiplies this error by $\beta-\alpha$ and gives  $\mathfrak{c}=2$ (see Lemma~\ref{lem:ellipseAB}). Therefore we obtain the factor $(\beta-\alpha)/2^{2n}= 4^{-n}\cdot (\beta-\alpha)$.
\qed
\end{proof}

\section{Petras' algorithm}\label{sec:Petras-algorithm}

In this section we recall an algorithm from \cite{P02} for the verified integration of piecewise analytic functions. We will call it \emph{the Petras' algorithm}.
There are three enigmatic points in the algorithm: (1) calculating an upper bound for $||f||_{\infty}$, (2) checking the analyticity and (3) boundedness of $f$. We assume that there exist subroutines (or oracles):
\begin{itemize}
\item {\tt UpperBound}$(f, a, b)$ which returns an upper bound $M$ for $||f||_{\infty} = \sup_{x\in [a,b]} |f(x)|$,

\item {\tt IsAnalytic}$(f,\alpha,\beta, A)$ such that
if it is true, then $f$ is analytic on $\ro(\alpha,\beta,A,\sqrt{A^2 - 1})$; we do not require that the converse is true,

\item {\tt ComplexBound}$(f, \alpha, \beta, A)$ returning a value greater than or equal to  $m(f;\alpha,\beta, A, \sqrt{A^2 - 1})$ provided that $f$ is analytic on $\ro(\alpha, \beta,A,\sqrt{A^2 - 1})$.

\end{itemize}

 In further analysis of Petras' algorithm we will formulate some conditions regarding the properties of the above functions, however we will not include their cost
 in the complexity estimations, even though they might be hard to compute and depend substantially on $f$. These issues will require the precise definition
 of computable analytic functions (see for example \cite{GaertnerHotz12, Hoeven05} and the references given there), which we do not consider in this paper.

\subsection{Formulation of the algorithm}

The input of the algorithm consists of:
\begin{itemize}
\item the integrand  $f$ and the interval of integration $[a, b]$; we require that $f$ is piecewise analytic and bounded on $[a,b]$,

\item an accuracy bound, $\varepsilon > 0$.
\end{itemize}
The algorithm also uses configuration constants:
\begin{itemize}
\item $A>1$ is used in the definition of the area $\ro$ (see Definition~\ref{def:rho-rectangle}) which is needed to compute {\tt IsAnalytic} and {\tt ComplexBound} functions (recall that $B = \sqrt{A^2 - 1}$),

\item the number $c>1$ for the estimation of function values.
\end{itemize}


The algorithm is as follows.
\begin{itemize}
\item[{\em In:}] $f, [a, b]$, $\varepsilon$

\item[1.] We choose a sequence $({\cal Q}_n)_{n\in \nat}$ of quadrature formulas to be Gaussian or Clenshaw-Curtis, i.e.\ we determine constants $D$ and $E$ (compare with~(\ref{eq:gauss-error}) and~(\ref{eq:c-c-error}) in Remark~\ref{rem:gauss-cc-errors}). These quadratures satisfy error estimates of the form
$$R_{n,[\alpha, \beta]}[f] \leqslant D\cdot E^{-n}\cdot (\beta - \alpha)m(f; \alpha, \beta)$$
for functions analytic on $\ro(\alpha,\beta,A,\sqrt{A^2 - 1})$.

\item[2.] $M$ := {\tt UpperBound}$(f, a, b)$;

\item[3.]
Assume that we have already partitioned $[a,b]$ into $k$ intervals, i.e.
$$[a,b]:= [a_0, a_1]\cup [a_1, a_2]\cup\ldots\cup [a_{k-1}, a_k]$$
and there is $J\subset \{1, \ldots, k\}$ such that
\begin{eqnarray}
\mbox{{\tt IsAnalytic}$(f, a_{j-1}, a_j,A)$}\mbox{ and \tt ComplexBound}(f, a_{j-1}, a_j, A) \leqslant cM, & & \mbox{for } j\in J\label{eqn:j-intervals}\\[1ex]
\sum_{j\not\in J}\left| a_j - a_{j-1}\right| > \frac{\varepsilon}{2M}.\label{eq:repeat-condition}
\end{eqnarray}
We choose among intervals not belonging to $J$ the longest interval
 and bisect it. Repeat step 3 as long as the condition~(\ref{eq:repeat-condition}) holds.

\item[4.] Now $\sum_{j\not\in J}\left| a_j - a_{j-1}\right| \leqslant \frac{\varepsilon}{2M}$ (condition~(\ref{eq:repeat-condition}) does not hold) and
\begin{enumerate}
\item[a)] for the intervals $[a_{j-1}, a_j]$ where $j\in J$ we calculate the integral using a quadrature formula  $Q_n$, where
\begin{equation}
n \geqslant \frac{1}{\ln E}\ln \frac{2D(b-a)cM}{\varepsilon}, \label{eq:number-of-points}
\end{equation}

\item[b)] for all remaining intervals we take the value of the integral to be 0.
\end{enumerate}

\item[{\em Out:}] $q := \sum_{j\in J} {\cal Q}_{n;[a_{j-1}, a_j]}[f]$
\end{itemize}


\begin{remark}\rm
The requirement in~(\ref{eqn:j-intervals}) that $f$ is analytic on $\ro(a_{j-1}, a_j,A,\sqrt{A^2 - 1})$ was missing in Petras' original paper~\cite{P02}.
\end{remark}


\subsection{Numerical accuracy}

\begin{lemma} \label{thm:error}
Assume that the Petras' algorithm stops returning $q$. Then
\begin{equation}
\left|\int_a^b f(x)dx - q \right| < \varepsilon.
\end{equation}
\end{lemma}

\begin{proof}
For each $[\alpha,\beta]$ for which $m(f;\alpha, \beta)\leqslant cM$ we want the error of the quadrature to be less than
$$(\beta - \alpha)\frac{\varepsilon}{2(b-a)},$$
so that summing the error over all intervals we obtain the global error on $[a,b]$ less than or equal to $\eps/2$.
We can calculate $n$, for all intervals, using an estimate
\begin{eqnarray*}
R_{n,[\alpha,\beta]}[f]  \leqslant  D\cdot E^{-n}\cdot (\beta - \alpha)m(f; \alpha,\beta)  \leqslant   D\cdot E^{-n}\cdot (\beta - \alpha)cM
   <  (\beta - \alpha)\frac{\varepsilon}{2(b-a)}
\end{eqnarray*}
and obtain
\begin{eqnarray*}
n & \geqslant & \frac{1}{\ln E}\ln \frac{2D(b-a)cM}{\varepsilon} .
\end{eqnarray*}

Now we can estimate the total error of the quadrature on all the intervals belonging to $J$:
\begin{eqnarray*}
\sum_{j\in J} R_{n;[a_{j-1}, a_j]}[f]  &\leqslant &  \sum_{j\in J}(a_j - a_{j-1})\frac{\varepsilon}{2(b - a)}
   \leqslant  \frac{\varepsilon}{2(b - a)}\sum_{j\in J}(a_j - a_{j-1})\\
  & \leqslant & \frac{\varepsilon}{2(b - a)} (b-a)   =  \frac{\varepsilon}{2}
\end{eqnarray*}

The error of the integration on the intervals not belonging to $J$ is not greater than
$M\cdot\varepsilon/(2M) = \varepsilon/2$ (an upper bound for $||f||_{\infty}$ times the total length of intervals where we do not use a quadrature).

Therefore we obtain
\begin{eqnarray*}
R_{n,[a,b]}[f] & = & \sum_{j\in J} R_{n;[a_{j-1}, a_j]}[f]  + \sum_{j\not\in J} R_{n;[a_{j-1}, a_j]}[f]   \leqslant  \frac{\varepsilon}{2} + \frac{\varepsilon}{2}
   =  \varepsilon.
\end{eqnarray*}
\qed
\end{proof}

\subsection{The algorithm stops under reasonable assumptions}

It is not obvious whether the algorithm stops. If the function {\tt ComplexBound} returns values which do not ``match'' the bound $c\cdot \sup_{x \in [a,b]} |f(x)|$ (for example, it always returns $2$, while $c\cdot \sup_{x \in [a,b]} |f(x)| \approx 1$) the algorithm might run forever.  Thus we have to assure that such case will never happen.
For this we need to provide certain assumptions regarding the quality (but not complexity) of functions {\tt IsAnalytic} and {\tt ComplexBound}.
\begin{definition}\label{def:compatibility}
Let $A>1$. We say that the algorithms {\tt IsAnalytic} and {\tt ComplexBound} satisfy the {\em compatibility condition} if for any $z \in [a,b]$ such that $f$ is analytic in some neighbourhood of $z$,
there exists an open set $U \subset \complex$, such that $z \in U$ and
\begin{equation*}
\mbox{if } [\alpha,\beta] \subset U \mbox{ then {\tt IsAnalytic}$(f,\alpha,\beta,A)$ and {\tt ComplexBound}$(f,\alpha,\beta,A) \leqslant c\cdot {\tt UpperBound}(f, a, b)$.}
\end{equation*}
\end{definition}

\begin{theorem}
Assume that $f$ is piecewise analytic and the Lebesgue measure of points in $[a,b]$ in which $f$ is not analytic is equal to zero. If functions {\tt IsAnalytic} and {\tt ComplexBound} satisfy the  compatibility condition, then Petras' algorithm stops.
\end{theorem}
\begin{proof}
From the assumptions about {\tt IsAnalytic} and {\tt ComplexBound} functions it follows that any interval created during the Petras' algorithm containing only the points of analyticity will be subdivided (possibly after several subdivisions) into smaller intervals that will eventually be accepted.

Hence the total length of bad intervals (not in $J$) goes to zero and the algorithm stops.
\qed
\end{proof}

\begin{remark}
From now on, suprema on the real axis, analyticity check and the functional $m$ should be understood as  values returned by the above subroutines. For the sake of simplicity we use $M$ to denote {\tt UpperBound}$(f, a, b)$.
\end{remark}

\subsection{The case of analytic functions}

\begin{lemma}\label{lm:analytic-fun}
Assume that $f$ is analytic on $[\alpha, \beta] \subset [a,b]$  and algorithms {\tt IsAnalytic} and {\tt ComplexBound} satisfy the compatibility condition. Then
the number of intervals accepted in the Petras' algorithm and covering $[\alpha,\beta]$ is finite and does not depend on $\eps$.
\end{lemma}
\begin{proof} It follows from the compactness of $[\alpha,\beta]$ and the compatibility condition (see Definition~\ref{def:compatibility}) that there exists $\delta$, such that any interval of length less than or equal to $\delta$ is accepted
by the Petras' algorithm.
\qed
\end{proof} 

\section{Tools for the analysis of the algorithm}
\label{sec:alg-analysis}
The Petras' algorithm has two parts: geometrical and computational. During the geometrical part (step~3) a partition of an interval $[a, b]$ is constructed based on a region of analyticity and boundedness~$\D$.
In the computational part (step~4), the integral is computed using a chosen quadrature on each of the resulting sub-intervals.

For a piecewise analytic function $f$ and $a \in \reals_+$ we define the region of analyticity and boundedness as
\begin{equation*}
  \D(f,a)= \overline{\left\{z \in \complex \ | \ z \in \doa{f}\wedge |f(z)| \leqslant a \right\}}.
\end{equation*}
Notice that in the perfect world, in the geometrical part the actual object of concern should be the set $\D = \D(f, c\cdot \sup_{x\in [a,b]}|f(x)|)$.
However, the set $\D$ is not known explicitly. Instead, the algorithm implicitly analyses the set $\D(f, cM)\subset \D$. It follows that (in our analysis) it is the geometric properties of $\D$ (and not the details of the function~$f$) that influence the actual complexity of the algorithm.

Further in the complexity estimations we use the following notations: intervals to be bisected are called ``bad'' and those that need not be bisected are ``proper."

The number of proper intervals generated by Petras' algorithm is denoted by
$$\mathcal{Z}(\D(f,c M), \eps).$$

Recall that on each of the proper intervals the quadrature (Gauss-Legendre or Clenshaw-Curtis) with

$$n = \left\lceil \frac{1}{\ln E}\ln \left(2Dc \cdot \frac{M(b-a)}{\varepsilon}\right) \right\rceil= \Theta(|\ln \varepsilon|)$$

points (see (\ref{eq:number-of-points}) in step 4 of the algorithm)  is calculated, thus the algorithm performs

\begin{equation*}
N(f,\eps)
   = \Theta(|\ln \varepsilon|) \cdot \mathcal{Z}(\D(f,c M), \eps)
\end{equation*}
evaluations of $f$. To complete the estimation of complexity we need to count the number of proper intervals generated by the algorithm.

For any $p >1$ we give examples of functions for which
\begin{equation*}
  \mathcal{Z}(\D(f,c M), \eps) = \Theta \left( \frac{1}{\eps^{p-1}} \right),
\end{equation*}
while the analysis in Petras' paper deals with the classes of functions for which
\begin{equation*}
   \mathcal{Z}(\D(f,c M), \eps) = \bigo (|\ln \varepsilon|).
\end{equation*}

\subsection{Petras-type conditions}
\label{subsec:petras-cond}
Let us fix $A$ and $B = \sqrt{A^2 - 1}$. In \cite{P02} Petras proposed a condition (we added the analyticity requirement missing in \cite{P02})
\begin{equation*}
\exists \gamma > 0\ \forall \alpha, \beta:\
\left[\ro(\alpha, \beta)\subset {\cal V}_{\gamma,S} \Longrightarrow  f \mbox{ is analytic on } \ro(\alpha,\beta) \mbox{ and }
 m(f; \alpha, \beta)\leqslant c \sup_{x \in [-1,1]} |f(x)|\right]
\end{equation*}
to estimate the complexity of the algorithm. We find this requirement too strong as it excludes a lot of functions.

Before presenting our generalization of the Petras' condition we state several theorems for the function $f(x)=\sin(1/x)$ for $x \in [-1,1]$, in order to illustrate that this condition does not hold and to motivate the use of the sets $D^p_{\gamma,S}$ in further analysis.

\begin{theorem}\label{thm:petras-insufficient} Let us take an arbitrary $Z>0$.   Let $f(z) = \sin(1/z)$. In the neighborhood of~$0$ the condition
\begin{equation}\label{eqn:petras-mod}
\exists \gamma > 0\ \forall \alpha, \beta\in [-1,1]:\
\left(\ro(\alpha, \beta)\subset {\cal V}_{\gamma,\{0\}} \quad\Longrightarrow\quad \sup_{z\in\ro(\alpha, \beta)}\left|f(z)\right|\leqslant Z\right)
\end{equation}
does not hold for any  $\gamma>0$,  i.e.
$$\forall \gamma > 0\ \exists \alpha, \beta\in [-1,1]:\ \left( \ro(\alpha, \beta)\subset {\cal V}_{\gamma,\{0\}} \wedge\ \sup_{z\in\ro(\alpha, \beta)}\left|f(z)\right|> Z\right).$$
\end{theorem}

\begin{proof}
Assume $z=x+iy$. Let $r^2=x^2 + y^2$. Then
\begin{equation}\label{eqn:1-z}
\frac{1}{z} = \frac{x-iy}{r^2}.
\end{equation}
Since
\begin{equation*}
  \sin(1/z)=\frac{\exp\left(i/z\right) - \exp\left(-i/z\right)}{2i}
\end{equation*}
we obtain
\begin{eqnarray*}
  \sin(1/z)
  & = & \frac{1}{2i}\left(
  \exp\left(i\frac{ x}{r^2}\right) \exp\left(\frac{y}{r^2}\right) - \exp\left(-i\frac{ x}{r^2}\right) \exp\left(-\frac{y}{r^2}\right)\right)
\end{eqnarray*}
and
\begin{equation}\label{eqn:abs}
|\sin(1/z)| \geqslant \frac{1}{2}\left(\exp\left(\frac{|y|}{r^2}\right) - \exp\left(-\frac{|y|}{r^2}\right)\right).
\end{equation}
For any $\gamma>0$ consider a point $z=x+i \gamma x$ with $x>0$. Since $y=\gamma x > 0$ there is
\begin{eqnarray*}
  \frac{|y|}{r^2}= \frac{\gamma x}{\gamma^2 x^2 + x^2}= \frac{\gamma}{x (\gamma^2 + 1)}
\end{eqnarray*}
and finally by~(\ref{eqn:abs}) we have
\begin{eqnarray*}
|\sin(1/z)| & \geqslant &  \frac{1}{2}\left(\exp\left(  \frac{\gamma}{x (\gamma^2 + 1)} \right) - \exp\left( -\frac{\gamma}{x (\gamma^2 + 1)} \right)\right)\\
    & \to & \infty , \quad \mbox{as } x \to 0.
\end{eqnarray*}
\qed
\end{proof}

Condition (\ref{eqn:petras-mod}) in the above theorem is not satisfied since on the border of ${\cal V}_\gamma$ the absolute value of $\sin (1/z)$ is too large. This is due to the fact that the term $\exp\left(|y|/r^2\right)$ grows to infinity as $z \to 0$. A workaround is to restrict the values of $\exp\left(|y|/r^2\right)$; however this leads to a different region, i.e.\ $D^2_{\gamma,S}$.

\begin{theorem}\label{thm:petras-condition-1}
Let $S = \{0\}$ and let $f(z) =  \sin (1/z)$.
Then for any $c > 1$
there exists $\gamma>0$ such that
\begin{equation}
  \sup_{z \in D^2_{\gamma} } |f(z)| \leqslant c \cdot \sup_{x \in [-1,1]}|f(x)|.
\end{equation}

Therefore  for any  $\alpha, \beta\in [-1,1]$
\begin{equation}\label{eqn:petras-condition-1}
\ro(\alpha, \beta)\subset \D^2_{\gamma} \quad\Longrightarrow\quad f \mbox{ is analytic on } \ro(\alpha, \beta)\ \wedge\ \sup_{z\in\ro(\alpha, \beta)}\left|f(z)\right|\leqslant c \cdot \sup_{x \in [-1,1]} |f(x)|.\\
\end{equation}

\end{theorem}

\begin{proof}
If $\ro(\alpha, \beta)\subset \D^2_{\gamma}$ then it is obvious that $f$ is analytic on $\ro(\alpha, \beta)$. Thus it is enough to check the second part of the thesis.

Assume $z=x+iy$. Let  $r^2 = x^2 + y^2$.
As in the proof of Theorem~\ref{thm:petras-insufficient} we have
\begin{eqnarray*}
  \sin(1/z) & = & \frac{1}{2i}
  \left(\exp\left(i\frac{ x}{r^2}\right) \exp\left(\frac{y}{r^2}\right) - \exp\left(-i\frac{ x}{r^2}\right) \exp\left(-\frac{y}{r^2}\right)\right)
\end{eqnarray*}
and then
\begin{eqnarray*}
  |\sin(1/z) | \leqslant \frac{1}{2} \left(\exp\left(\frac{|y|}{r^2}\right) + \exp\left(-\frac{|y|}{r^2}\right)\right).
\end{eqnarray*}

For $\gamma > 0$ let us define $\W_\gamma$ by
\begin{equation*}\label{eqn:disk}
  \W_\gamma  = \left\{ x+iy: \  \frac{|y|}{r^2} \leqslant \gamma \right\}. 
\end{equation*}
We have (because function $x \mapsto x + \frac{1}{x}$ is increasing for $x >1$)
\begin{equation*}
   |\sin(1/z) | \leqslant \frac{1}{2} \left(\exp(\gamma) + \exp(-\gamma) \right), \quad z \in \W_\gamma.
\end{equation*}
Since
\begin{equation*}
\frac{1}{2}\lim_{\gamma\to 0} \left(\exp(\gamma) + \exp(-\gamma) \right) = 1,
\end{equation*}
by taking $\gamma$ sufficiently close to 0 we get
\begin{equation*}
\sup_{z \in \W_\gamma} \left|  \sin \left(1/z \right)\right|
< c \cdot \sup_{x \in [-1,1]} | \sin(1/x)|.
\end{equation*}

To complete the proof it is enough to show that
 $ \D^2_{\gamma} \subset  \W_\gamma.$
Let $y\geqslant 0$ (the other case is symmetric). It is easy to see that $\W_\gamma$ describes
the complement of a disk of radius $\frac{1}{2 \gamma}$ centred at $\left(0,\frac{1}{2 \gamma}\right)$, i.e.\
\begin{equation*}
z=x+iy  \in \W_\gamma \quad\Longleftrightarrow\quad  \left(y - \frac{1}{2 \gamma}\right)^2  + x^2 \geqslant \frac{1}{4 \gamma^2}.
\end{equation*}
Observe that $\W_\gamma$ contains the following set
\begin{equation*}
\widetilde{\W}_\gamma= \left\{ x+iy: \  |x| \leqslant \frac{1}{2 \gamma}; \quad |y| \leqslant \frac{1- \sqrt{1-4 \gamma^2 x^2}}{2 \gamma}\right\}  \cup \left\{x+iy: \ |x| > \frac{1}{2\gamma}\right\}
\end{equation*}
and notice that for $|x| \leqslant \frac{1}{2 \gamma}$
$$
\frac{1- \sqrt{1-4 \gamma^2 x^2}}{2 \gamma} >  \gamma x^2.
$$
This shows that $\D^2_{\gamma} \subset \widetilde{\W}_\gamma \subset  \W_\gamma $.
\qed
\end{proof}

 Theorem~\ref{thm:petras-condition-1} says that for $z\mapsto  \sin(1/z)$ there exists a region $\D^2_{\gamma}$ where the function is analytic and appropriately bounded. The next theorem says the opposite:
there exists a region of the same shape as before, but such that (on the boundary of this region
close to the singular point) the values of the function are arbitrarily large.

\begin{theorem}\label{thm:petras-condition-2}
Let $S = \{0\}$ and let $f(z) = \sin(1/z)$.
For any $c > 1$ there exists $\gamma > 0$ such that for any $\alpha\in [-1/\gamma,1/\gamma]$ and $\beta\in (\alpha, 1]$
\begin{equation}\label{eqn:petras-condition-2}
\ro(\alpha, \beta)\not\subset \D^2_{\gamma} \quad\Longrightarrow\quad f \mbox{ is not analytic on } \ro(\alpha, \beta)\ \vee \sup_{z\in\ro(\alpha, \beta)}\left|f(z)\right| > c.
\end{equation}
\end{theorem}

\begin{proof}
Let us fix $c>1$. If $\alpha < 0 < \beta$ then $f$ is not analytic on $\ro(\alpha, \beta)$.

Let us consider the case when $f$ is analytic on $\ro(\alpha, \beta)$. Notice that if $\ro(\alpha, \beta)\not\subset \D^2_{\gamma}$ then $\ro(\alpha, \beta) \cap \partial \D^2_{\gamma} \neq \emptyset$,  where $\partial \D^2_{\gamma}$ is a border of $\D^2_{\gamma}$. Thus let us consider the region $${\cal W}_\gamma = \partial \D^2_{\gamma} \cap \{z:\ {\rm Re}(z) \in [-1/\gamma,1/\gamma]\} \setminus \{0\}.$$
Then
$$\sup_{z\in\ro(\alpha, \beta)}\left|\sin(1/z)\right| \geqslant
\inf_{z\in {\cal W}_\gamma }\left|\sin(1/z)\right|.$$

Assume $z = x + iy\in {\cal W}_{\gamma}$ and let $r^2 = x^2 + y^2$. Then  (for $|x|\leqslant 1/\gamma$)
\begin{equation*}
\frac{|y|}{r^2} = \frac{\gamma x^2}{x^2 + (\gamma x^2)^2} = \frac{\gamma}{1 + \gamma^2 x^2} \geqslant \frac{\gamma}{2}.
\end{equation*}
Since (see proof of Theorem~\ref{thm:petras-insufficient})
\begin{eqnarray*}
\left|\sin(1/z)\right| & \geqslant & \frac{1}{2}\left(\exp\left( \frac{|y|}{r^2}\right) - \exp\left(- \frac{|y|}{r^2}\right) \right)
\end{eqnarray*}
and $x\mapsto x - \frac{1}{x}$ is an increasing function we have
\begin{equation*}
\inf_{z\in {\cal W}_{\gamma}}\left|\sin(1/z)\right|\geqslant
\frac{1}{2} \left(\exp\left( \frac{\gamma}{2}\right) - \exp\left(- \frac{\gamma}{2}\right) \right).
\end{equation*}
Now, for any $c>1$ it is enough to take $\gamma > 0$ such that
$\inf_{z\in {\cal W}_\gamma} |\sin(1/z)| > c$, i.e.
\begin{eqnarray*}
\gamma & > & 2 \ln (c + \sqrt{1 + c^2})
\end{eqnarray*}
and the condition~(\ref{eqn:petras-condition-2}) holds.
\qed
\end{proof}

The above investigations justify why we consider regions bounded by curves $\gamma x^p$.
We do realize that this choice is arbitrary and non-exhaustive.
Nevertheless it allows us to show that region of analyticity is an important parameter in the complexity investigations.  In Section~\ref{sec:examples} we consider  the functions $f(z) = z^k\sin(1/z)$ which need regions of the form  $D^p_{\gamma}$.

Although so far we have considered only $p \geqslant 1$, we can in fact state definitions for any $p\in\reals_+$.
We present modified Petras' condition introducing the notion of order and a converse implication. Notice that in Definitions~\ref{def:petras-routines} and~\ref{def:petras-negative}, $\ro$ and $m$ are implicitly dependent on $A$ and $B$.

\begin{definition}
\label{def:petras-routines}
Let us fix  $A$, $B$ and $c>1$.
We say that \emph{a function $f$ satisfies positive Petras' condition of order $p$} (abbreviated as \PPC{p, \gamma, S}) on $[a,b]$, if there exist $S = \{s_1, \ldots, s_m\}$ with
$a \leqslant s_1 < \ldots < s_m \leqslant b$ and $\gamma > 0$ such that
\begin{align}
\ro(x, y)\subset \D^p_{\gamma,S} \Longrightarrow f \mbox{ is analytic on } \mathbf{\ro}(x, y) \mbox{ and } m(f; x, y)\leqslant c M,\quad\quad\mbox{}\label{eq:Petras-assumption}\\
\mbox{ for all $x$, $y$ such that $[x,y] \subset [a,b]$.}\nonumber
\end{align}
\end{definition}

\begin{definition}
\label{def:petras-negative}
Let us fix  $A$, $B$ and $c>1$.
We say that \emph{a function $f$ satisfies negative Petras' condition of order $p$} (\NPC{p, \gamma, s, \beta}) on $[a,b]$, if
there exist $s\in [a,b]$, $\gamma > 0$ and $\beta >0$ such that one of the following two conditions is satisfied
\begin{align}
\ro(x, y)\not\subset \D^p_{\gamma} \Longrightarrow
 \mbox{ either } (f \mbox{ is not analytic on } \ro(x, y)) \mbox{ or } m(f; x, y)> c M,\quad\quad\mbox{}
\label{eq:Petras-assumption-opp}\\
\mbox{for all $x$, $y$ such that $x\in [s, s+\beta]$ and $y\in (x, b]$,}\nonumber
\end{align}
\begin{align}
\ro(x, y)\not\subset \D^p_{\gamma} \Longrightarrow
 \mbox{ either } (f \mbox{ is not analytic on } \ro(x, y)) \mbox{ or } m(f; x, y)> c M,\quad\quad\mbox{}
\label{eq:Petras-assumption-opp-left}\\
\mbox{for all $x$, $y$ such that  $y\in (s-\beta,s]$ and $x\in [a,y]$,}\nonumber
\end{align}
We refer to (\ref{eq:Petras-assumption-opp}) as  (\NPC{p, \gamma, s, \beta})-right and (\ref{eq:Petras-assumption-opp-left}) as    (\NPC{p, \gamma, s, \beta})-left.
\end{definition}
Notice that in the above definition $\D^p_{\gamma} = \D^p_{\gamma, \{s\}}$.


The idea (how we use PPC and NPC conditions to estimate the number of proper intervals) is presented in Remark~\ref{ex:idea} and a more detailed treatment follows in Section~\ref{subsec:Various-regions}.
Notice that, in general, we want to point out the classes of functions for which the complexity of Petras' algorithm is $\Theta(|\ln\eps|/\eps^{p-1})$. Therefore we consider only those functions for which NPC and PPC are satisfied. To show that the complexity of the algorithm cannot be better it is enough to show one singular point of $f$ where the number of proper intervals generated by the algorithm is $\Omega(1/\eps^{p-1})$ --- a lower bound. To show that the complexity of the algorithm is not worse we have to prove that in any point of $S$ the algorithm cannot produce more proper intervals than $\bigo(1/\eps^{p-1})$ --- an upper bound.

\begin{remark}\label{ex:idea}
Let $\D = \D(f,c M)$. Assume that a function $f$ satisfies \PPCS{p, \gamma, \{s\}} and \NPC{p', \gamma', s, \beta} (with $\gamma' \geqslant \gamma$, $p' \leqslant p$) on the whole interval, i.e.\ $x\in [a,b]$ in~(\ref{eq:Petras-assumption-opp},\ref{eq:Petras-assumption-opp-left}), therefore $\D_{\gamma}^p\subset \D\subset \D_{\gamma'}^{p'}$.
If we consider intervals $[\alpha, x]$, $[\alpha, \alpha']$ and $[\alpha, y]$ such that $\ro(\alpha,x) \subset \D_{\gamma}^p$, $\ro(\alpha, \alpha') \subset \D$ and
$(\ro(\alpha, y)\not\subset {\rm int} \D_{\gamma'}^{p'} \wedge \ro(\alpha, y)\subset \D_{\gamma'}^{p'})$,
then $[\alpha,x]\subset [\alpha, \alpha']\subset [\alpha, y]$.
Thus investigating intervals in $\D_\gamma$ (they are the shortest and therefore their number is the largest) we are able to estimate the number of proper intervals from above. And while investigating intervals on the edge of $\D_{\gamma'}$ (the longest ones) we are able to estimate the number of proper intervals from below.

\end{remark}

\subsection{Geometric lemmas}\label{subsec:geo-gen}
Let us fix $S=\{0\}$. We define
\begin{equation*}
\D^{\Gamma} = \{x + i y:\ |y|\leqslant \Gamma(x)\},
\end{equation*}
where $\Gamma:\ [0, \infty)\to [0, \infty)$ is continuous, $\Gamma(0) = 0$ and $\Gamma$ is strictly increasing for $x > 0$.
The goal of this section is to find, for a fixed $x$ or $y$,  the largest possible $d=y-x$ such that $\ro(x,y,A,B) \subset \D^\Gamma$.

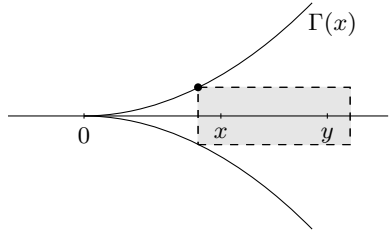
\begin{figure}[htb]
\centering
\begin{tikzpicture}
\filldraw[fill=gray!20!white, draw=gray!20!white ] (2.5,-0.38) rectangle (4.5,0.38);
\draw (1,0) parabola (4.0,1.5) node[anchor=north] {\small \hspace*{4ex}$\Gamma(x)$};
\draw (1,0) parabola (4.0,-1.5);
\draw (1 cm,1pt) -- (1 cm,-1pt) node[anchor=north] {\small $0$};

\draw (0,0) -- (5,0);
\draw[dashed] (2.5,-0.38) rectangle (4.5,0.38);
\draw[dashed] (2.5,-0.38) rectangle (4.5,0.38);
\draw (2.8,1pt) -- (2.8,-1pt) node[anchor=north] {\small $x$};
\draw (4.2,1pt) -- (4.2,-1pt) node[anchor=north] {\small $y$};
 \fill (2.5, 0.38) circle [radius=1.5pt];
\end{tikzpicture}
\caption{Points $x$ and $y$ are chosen so that the top left corner of $\ro(x, y, A, B)$ belongs to the line $x\mapsto \Gamma(x)$.}\label{fig:situation}
\end{figure}

By the definition of $\ro$ the top left corner of $\ro(x, x+d, A, B)$ is at
$$\left(x - \frac{d}{2}(A - 1), B\frac{d}{2}\right)$$
and we want this point to lie on the line $x\mapsto \Gamma(x)$, thus (see Figure~\ref{fig:situation})
\begin{equation}\label{eq:G-d-left}
\Gamma\left(x - \frac{d}{2}(A - 1)\right) = B\frac{d}{2}.
\end{equation}
This is the desired formula for $d=y-x$ parametrized by $x$, the left end of the interval.
To obtain a formula parametrized by $y$ we substitute $y - d$ for $x$ in equation~(\ref{eq:G-d-left}) and obtain
\begin{equation}\label{eq:G-d-right}
\Gamma\left(y - \frac{d}{2}(A + 1)\right) = B \frac{d}{2}.
\end{equation}

\begin{theorem}
\label{thm:G-d-exists}
There exist strictly increasing, continuous functions $d_L, d_R: [0,\infty) \to [0,\infty)$ solving equations  (\ref{eq:G-d-left}), (\ref{eq:G-d-right}), respectively.
Moreover:
\begin{eqnarray}
   d_L(x) & < & \frac{2}{B} \Gamma(x), \quad x >0,   \label{eq:deltaL-Upper}  \\
    d_R(y) & < & \frac{2}{B} \Gamma(y), \quad y >0.  \label{eq:deltaR-Upper}
  \end{eqnarray}
\end{theorem}
\begin{proof}
It is clear that $d_L$ is strictly increasing in $x$. It is enough to translate  $\ro(x,y,A,B)$, with top left corner on the line $x\mapsto \Gamma(x)$, to the right. The shifted rectangle will be contained in the interior
of $\D^\Gamma$. Similarly, $d_R$ is increasing in $y$.

For the proofs of (\ref{eq:deltaL-Upper}) and (\ref{eq:deltaR-Upper}) observe that for $x>0$ it holds
\begin{eqnarray*}
  \frac{B}{2} d_L(x) & = &\Gamma \left(x - \frac{A-1}{2}d_L(x) \right) < \Gamma \left(x  \right)\\
   \frac{B}{2} d_R(y) & = & \Gamma \left(y - \frac{A+1}{2}d_R(x) \right) < \Gamma (y).
\end{eqnarray*}
\qed
\end{proof}

We would like to obtain bounds of the following form
\begin{equation}
 c_2 \Gamma(x) \leqslant d_L(x), d_R(x) \leqslant c_1 \Gamma(x),
\end{equation}
for some $c_1, c_2 > 0$. The existence of $c_1$ follows from Theorem~\ref{thm:G-d-exists}.
The existence of $c_2$ is treated in Section~\ref{subsec:regions-p-bigger-1} for $\Gamma(x)=\gamma x^p$, $p>1$.
It turns out that for $p\leqslant 1$ this is not true; in fact we obtain linear estimates from above and from below (see Theorem~\ref{thm:d(x)estm-p<1} and equation~(\ref{eq:d-p=1})).

For the case $\Gamma(x) = \gamma x^p$ let us set
\begin{equation*}
  h=\frac{B}{2\gamma}.
\end{equation*}
Then observe that equations (\ref{eq:G-d-left}) and (\ref{eq:G-d-right}) have the following form
\begin{equation}\label{eqn:next-point1}
(x - g\cdot d)^p = h\cdot d,
\end{equation}
provided that we set either
\begin{equation*}
g = g_L = (A - 1)/2, 
\end{equation*}
or
\begin{equation*}
 g = g_R = (A+1)/2, 
\end{equation*}
  respectively.
Equation (\ref{eqn:next-point1}) defines implicitly a function $d(x)$.  In the following subsections we will estimate $d(x)$.

\subsubsection{The case $p > 1$}\label{subsec:regions-p-bigger-1}
Our goal in this section is to develop some estimates for the solution of (\ref{eqn:next-point1}) for $p>1$.
To develop intuitions consider an integer $p$ and
a series expansion of $d(x) = d_0 + d_1 x + d_2 x^2 + \ldots$\ .
Regrouping the terms in~(\ref{eqn:next-point1}) and taking $d_0 = 0$ we obtain
\begin{eqnarray*}
d(x)
    & = & \frac{1}{h}x^p + d_{p+1}x^{p+1} + ... \\
    & = & \frac{1}{h}x^p(1 + d_{p+1} x + ...)\\
    & = & \frac{c(x)}{h}\cdot x^p,
\end{eqnarray*}
where $c(x) = O(1)$ for small $x$.
These considerations lead us to a hypothesis that
\begin{equation}\label{eqn:d(x)-1}
d(x) = \frac{c(x)}{h}\cdot x^p
\end{equation}
for any $p>1$, such that there exist $c_1$, $c_2$ and
\begin{equation*}
0 < c_2 \leqslant c(x) \leqslant c_1,
\end{equation*}
for a bounded range of $x$.
Therefore from (\ref{eqn:next-point1}) and (\ref{eqn:d(x)-1}) we have the following implicit equation for $c(x)$
\begin{eqnarray*}
\left(1 - g\cdot \frac{c(x)}{h}\cdot x^{p-1}\right)^p & = & c(x).
\end{eqnarray*}

For given $g$, $h$ and $p > 1$ define
\begin{equation}\label{eqn:x-p>1}
\tilde{x} = \left(\frac{h}{g}\right)^{\frac{1}{p-1}}.
\end{equation}

\begin{lemma}\label{lm:bounds}
Let us consider the equation
\begin{equation}\label{eqn:to-solve}
F(x,c(x)) = \left(1 - c(x)\cdot \frac{g}{h}\cdot x^{p-1}\right)^p - c(x) = 0.
\end{equation}
This equation has exactly one solution $c(x) \in [0,1]$ which is continuous on $[0, \tilde{x}]$ and
$$\forall x\in [0, \tilde{x}]\ \exists c', c'': (0 < c'' < c' = 1) \wedge (c''\leqslant c(x) \leqslant c').$$
\end{lemma}

\begin{proof}
Notice that for every $x$,
$F(x,1) < 0$ and $F(x,0) = 1 > 0$ thus there exists $c\in (0,1]$ such that $F(x,c) = 0$.
We now show the uniqueness of $c$:
\begin{eqnarray*}
\frac{\partial F}{\partial c} & = & p\left(1 - c\cdot \frac{g}{h} x^{p-1} \right)^{p-1} \left(-\frac{g}{h} x^{p-1} \right) - 1 < 0.
\end{eqnarray*}
The function $c:[0, \tilde{x}]\to (0,1]$ is continuous. Since $[0, \tilde{x}]$ is compact, there exists $x_0\in [0, \tilde{x}]$ such that $$\forall x\in [0, \tilde{x}]:\ 0 < c(x_0) \leqslant c(x)$$
hence $c'' = c(x_0)$ is the bound we need.
\qed
\end{proof}

The following theorem gives us the desired upper and lower estimates for $d(x)$.
\begin{theorem}\label{thm:d(x)estm-p>1}
Let $p>1$ and $d(x)$ be the smallest nonnegative solution  of (\ref{eqn:next-point1}).
Then there exist $0 < c_2 < c_1$ such that for $x \in [0,1]$
\begin{equation}\label{eqn:d(x)-estm}
   c_2\cdot x^{p} \leqslant d(x) \leqslant c_1\cdot x^{p}.
\end{equation}
\end{theorem}
\begin{proof}
By~(\ref{eqn:d(x)-1}) we have $d(x) = x^p c(x)/h$.
Let  $\tilde{x}$ be as in~(\ref{eqn:x-p>1}), then by Lemma~\ref{lm:bounds} we immediately obtain:
\begin{equation*}
   \tilde{c}_2\cdot x^{p} \leqslant d(x) \leqslant \tilde{c}_1\cdot x^{p}, \quad x \in [0,\tilde{x}]
\end{equation*}
for some $0< \tilde{c}_2 < \tilde{c}_1$.
For $x\in [\tilde{x}, 1]$ by Theorem~\ref{thm:G-d-exists} we have:
\begin{itemize}
\item $d(x)$ is strictly increasing, thus there exists $\bar{c}_2 > 0$ such that $d(x) \geqslant \bar{c}_2 x^p$;
\item $d(x) < 2\gamma x^p/ B$, thus there exists $\bar{c}_1 > \bar{c}_2$ such that $d(x) \leqslant \bar{c}_1\cdot x^{p}$.
\end{itemize}
Taking $c_1 = \max \{\tilde{c}_1, \bar{c}_1\}$ and $c_2 = \min \{\tilde{c}_2, \bar{c}_2\}$ we obtain the thesis.
\qed
\end{proof}

\subsubsection{The case $p < 1$}\label{subsec:regions-p<1}

Consider the equation~(\ref{eqn:next-point1}) for $0 < p < 1$:
\begin{eqnarray}
(x - g\cdot d)^p & = & h\cdot d\nonumber\\
x - g\cdot d & = & h^{\frac{1}{p}}\cdot d^{\frac{1}{p}}\label{eqn:next-point2}.
\end{eqnarray}

Similar considerations as before for $p > 1$
lead us to a hypothesis that
\begin{equation}\label{eqn:d(x)-2}
d(x) = \frac{x}{g}\left(1 - c(x)\cdot \frac{h^{\frac{1}{p}}}{g^{\frac{1}{p}}}\cdot x^{\frac{1}{p} - 1}\right)
\end{equation}
where $c(x)$ is a bounded positive function.
Substituting~(\ref{eqn:d(x)-2}) for $d$ in~(\ref{eqn:next-point2}) we obtain the following implicit equation for $c(x)$
\begin{equation}\label{eqn-c(x)-2}
c(x) = \left(1 - c(x)\cdot \frac{h^{\frac{1}{p}}}{g^{\frac{1}{p}}}\cdot x^{{\frac{1}{p}}-1}\right)^{\frac{1}{p}}.
\end{equation}

As in the case $p > 1$ let us define, for given $g$, $h$ and $p < 1$
\begin{equation}\label{eqn:x-p<1}
\tilde{x} = \left(\frac{g}{h} \right)^{\frac{1}{1-p}}.
\end{equation}
Now, for $p < 1$ we have an analogue of Lemma~\ref{lm:bounds}.
\begin{lemma}\label{lm:bounds-2}
Let us consider the equation
\begin{equation}\label{eqn:to-solve-2}
F(x,c(x)) = \left(1 - c(x)\cdot \frac{h^{\frac{1}{p}}}{g^{\frac{1}{p}}}\cdot x^{{\frac{1}{p}}-1}\right)^{\frac{1}{p}} - c(x) = 0.
\end{equation}
This equation has exactly one solution $c(x)\in [0,1]$ which is continuous on $[0, \tilde{x}]$ and
$$\forall x\in [0, \tilde{x}]\ \exists c', c'': (0 < c'' < c' = 1) \wedge (c''\leqslant c(x) \leqslant c').$$
\end{lemma}

From (\ref{eqn:d(x)-2}) and the above lemma we obtain the following bounds for $d(x)$.
\begin{theorem}\label{thm:d(x)estm-p<1}
Let $p < 1$ and $d(x)$ be the smallest nonnegative solution of (\ref{eqn:next-point1}).
Then there exist $0 < c_2 < c_1$ such that for $x \in [0,\tilde{x}]$, where $\tilde{x}$ is as in~(\ref{eqn:x-p<1}), 
\begin{equation}\label{eqn-p-less-1}
x\left( \frac{1}{g} - c_1\cdot x^{\frac{1}{p}-1}\right) \leqslant d(x) \leqslant x\left( \frac{1}{g} - c_2\cdot x^{\frac{1}{p}-1}\right).
\end{equation}
\end{theorem}

\subsubsection{The case $p = 1$}

The equation~(\ref{eqn:next-point1}) for $p = 1$ has a form
\begin{eqnarray}
x - g\cdot d & = & h\cdot d\label{eqn:next-point3},
\end{eqnarray}
thus we obtain
\begin{equation}
d(x) = \frac{x}{g+h}. \label{eq:d-p=1}
\end{equation}

\begin{remark}
In further analysis we will refer to Theorem~\ref{thm:d(x)estm-p<1} with $p = 1$ (as it holds for $p\leqslant 1$). Thus it is sufficient to consider two cases, $p>1$ and $p\leqslant 1$.
\end{remark}

\section{Lower bound under NPC condition}
\label{sec:lower-bnd-estm}

We assume
\begin{itemize}
\item $[a,b] = [-1,1]$,
\item $\D = \D(f,c M)$,
\item $f$ satisfies \NPC{p', \gamma', s, \beta}-right (with $s\in S$, $\gamma' \geqslant \gamma$, $p' \leqslant p$).
\end{itemize}
Therefore $\D \cap \{ x+iy \ | \ x \in [s,s+\beta] \}\subset \D_{\gamma'}^{p'} \cap \{ x+iy \ | \ x \in [s,s+\beta] \}$. For the sake of simplicity of calculations we assume that $s = 0$.


\subsection{Rightward flow}\label{subsec:Various-regions}

We are interested in the number of proper intervals generated in the segment $[s,\beta]$ by the Petras' algorithm (computing the integral attaining global precision $\eps$ on $[a,b]$). In this section we define a rightward flow and we estimate this number from below by the number of steps in this flow. For the estimation from above we define, in Section~\ref{subsec:leftward-flow}, a leftward flow.

We will refer to the manner in which the algorithm creates successive intervals as a {\em Petras' process}.
During this process we obtain a sequence of proper intervals
\begin{equation}\label{petras-seq}
s\leqslant \alpha_0 < \alpha_1 < \ldots < \alpha_n, \quad\mbox{with}\quad \alpha_{n-1} < \beta \leqslant \alpha_n,
\end{equation}
where $\displaystyle \alpha_0 \leqslant \frac{\eps}{2M}$.\\[2ex]

The \NPC{p', \gamma', s, \beta}-right condition implies that there exist $s$ and $\beta$, such that \NPC{p', \gamma', s, \beta} holds for $x\in [s, s+\beta]$. Recall that we use $s = 0$ thus $\beta$ can be treated as a distance from the singular point.
\begin{enumerate}
\item We start from $x_0 = \alpha_0$, where $\alpha_0$ is
the point (close to $s$) where Petras' process has finished bisecting intervals.

\item Each new interval $[x_i, x_{i+1}]$ is created in such a way that $\ro(x_i, x_{i+1}) \subset \D^{p'}_{\gamma'}$ and upper left corner of $\ro(x_i, x_{i+1})$  is on the boundary of  $\D^{p'}_{\gamma'}$. Observe that all longer intervals $[x_i, x_{i+1}+\delta]$, $\delta>0$ lack this property.

\item We count the number of steps needed to exit the interval $[s,\beta]$ to the right, i.e.\ to have $x_{m-1} < \beta \leqslant x_m$.
\end{enumerate}

This process creates a sequence of intervals
\begin{equation}\label{right-seq}
\alpha_0  = x_0 < x_1 < \ldots < x_{m-1} < \beta \leqslant x_m.
\end{equation}

Each successive point $x_{n+1}$ is chosen as
\begin{equation}\label{eqn:process-general}
x_{n+1} = x_n + d_L(x_n),
\end{equation}
where function $d_L(x)$ is equal to $d(x)$ that solves equation (\ref{eqn:next-point1}) with
\begin{equation}
h=\frac{B}{2\gamma'}, \qquad g=g_L=\frac{A-1}{2}. \label{eq:gh-rightflow}
\end{equation}
By (\ref{eq:Petras-assumption-opp}) if $x_n \leqslant \beta$, then on any rectangle $\ro(w_1,w_2)$ containing $\ro(x_n,x_{n+1})$, either the function $f$ is not analytic, or $\sup_{z \in \ro(x_n,x_{n+1})} |f(z)| > cM$, so that $[w_1,w_2]$ will not be accepted by the Petras' algorithm as a proper interval.

Estimates for $d_L(x)$ are given by Theorem~\ref{thm:d(x)estm-p>1} (for $p>1$), Theorem~\ref{thm:d(x)estm-p<1} (for $p<1$) and (\ref{eq:d-p=1}) (for $p=1$).

\subsection{Estimation of the number of proper intervals from below}
Below we consider functions that satisfy \NPC{p, \gamma, 0, \beta}-right (Lemma~\ref{lm:lower-bound} and~\ref{lm:n-estimation-p<1}). The case of functions satisfying \NPC{p, \gamma, 0, \beta}-left is analogous.

\subsubsection{Case $p>1$}
\begin{lemma}\label{lm:lower-bound}
Assume that $p > 1$ and $f$ satisfies \NPC{p, \gamma, 0, \beta}-right. Then, for sufficiently small $\eps$, in the segment $[0, \beta]$
$$\mathcal{Z}(\D, \eps) = \Omega\left(\left(\frac{1}{\eps} \right)^{p-1} \right).$$
\end{lemma}

\begin{proof}
Let $\displaystyle\frac{\eps}{2M} < \beta$.
Let $t$ be the number of proper intervals in $\displaystyle\left[0,\beta\right]$.
From~(\ref{eqn:process-general}) and  Theorem~\ref{thm:d(x)estm-p>1} it follows that there
exists $c_1 > 0$, such that
\begin{equation*}
d_R(x) \leqslant c_1 \cdot x^p,
\end{equation*}
therefore
\begin{eqnarray*}
\alpha_{i+1} - \alpha_{i}  <  d_R(\alpha_i) \leqslant  c_1 \cdot \alpha_i^p.
\end{eqnarray*}
It is easy to see that if $x(t)$ is a solution of $x' = c_1 \cdot x^p$ with an initial condition
$\displaystyle x(0)=\frac{\eps}{2M}$, then
\begin{equation*}
  x(i) \geqslant \alpha_i, \quad  i=0,\dots, t.
\end{equation*}
The solution of $x' = c_1 \cdot x^p$ is given by
\begin{equation*}
  x(t)=\frac{x(0)}{\left[1 - x(0)^{p-1}c_1 (p-1)t \right]^{\frac{1}{p-1}}},
\end{equation*}
hence we can calculate the exit time to the right from
$\displaystyle \left[\frac{\eps}{2M}, \beta\right]$ (taking $x(t) = \beta$), i.e.
\begin{eqnarray*}
t & \geqslant & \frac{1}{\eps^{p-1}}\cdot\frac{(2M)^{p-1}}{c_1(p-1)} - \frac{1}{c_1(p-1)\beta^{p-1}}.
\end{eqnarray*}
Therefore (see Definition~\ref{def:asym-rate-groth})
\begin{equation}
\mathcal{Z}(\D, \eps) = \Omega\left(\left(\frac{1}{\eps} \right)^{p-1} \right).\label{lower-bound}
\end{equation}
\qed
\end{proof}

\subsubsection{Case $p\leqslant 1$}

\begin{lemma}\label{lm:n-estimation-p<1}
Assume that $p \leqslant 1$ and $f$ satisfies \NPC{p, \gamma, 0, \beta}-right. Then, for sufficiently small $\eps$, in the segment $[0, \eta]$, where $\eta = \min\{\tilde{x}, \beta\}$, for $\tilde{x}$ as in~(\ref{eqn:x-p<1}),
$$\mathcal{Z}(\D, \eps) = \Omega\left(\ln \frac{1}{\eps} \right).$$
\end{lemma}

\begin{proof}
Let $\eta = \min \{\tilde{x}, \beta\}$ and $\displaystyle \frac{\eps}{2M} < \eta$.
Let $k$ be the number of proper intervals in $[0,\eta]$.
Following (\ref{eq:gh-rightflow}) we set
\begin{equation} \label{eq:gr}
  g_R=\frac{A-1}{2}.
\end{equation}

By Theorem~\ref{thm:d(x)estm-p<1}  it follows that there exists  $c_2$ such that for $x\in [0, \tilde{x}]$
\begin{equation*}
d_R(x) \leqslant x\left(\frac{1}{g_R} - c_2 \cdot x^{\frac{1}{p}-1}\right) \leqslant \frac{x}{g_R}.
\end{equation*}

Now, for $\alpha_i\in [0, \eta] $, we have
\begin{equation*}
\alpha_{i+1} - \alpha_i  \leqslant  d_R(\alpha_i) \leqslant \frac{\alpha_i}{g_R}
\end{equation*}
thus
\begin{equation*}
\alpha_{i+1} \leqslant \alpha_i\left(1 + \frac{1}{g_R} \right)
\end{equation*}
and, since $\displaystyle\alpha_0\leqslant \frac{\eps}{2M}$, we obtain
\begin{equation*}
\alpha_i \leqslant \alpha_0\left(1 + \frac{1}{g_R}\right)^i\leqslant \frac{\eps}{2M}\left(1 + \frac{1}{g_R}\right)^i.
\end{equation*}
Therefore in particular
\begin{equation*}
\eta\leqslant \frac{\eps}{2M}\left(1 + \frac{1}{g_R}\right)^k
\end{equation*}
and
\begin{equation*}
  k \geqslant \frac{\ln \left(\frac{2 M\eta}{\eps}\right)}{\ln \left( 1+ \frac{1}{g_R}\right)} - 1,
\end{equation*}
thus according to the Definition~\ref{def:asym-rate-groth} we obtain the thesis.
\qed
\end{proof}

\section{Upper bound from PPC condition}
\label{sec:upper-bound}
We assume:
\begin{itemize}
\item $[a,b] = [-1,1]$,
\item $\D = \D(f,c M)$,
\item $f$ satisfies \PPCS{p, \gamma, S}, therefore $\D_{\gamma}^p\subset \D$.
\end{itemize}

\subsection{Modified Petras' algorithm}
 We consider a modified Petras' algorithm (abbreviated MPA) in which a decision if an interval is proper or bad is made using PPC condition instead of checking  {\tt IsAnalytic} and {\tt ComplexBound}. Additionally, if MPA decides that it is still running (the length of all bad intervals is greater than $\eps/(2 M)$), then it bisects all bad intervals in step 3 in the Petras' algorithm before checking condition (\ref{eq:repeat-condition}); thus in MPA all bad intervals have the same length.

Notice that any interval accepted by MPA is contained in an interval accepted by Petras' algorithm. Thus we can state:
\begin{enumerate}
\item a sum of lengths of proper intervals generated in Petras' algorithm is greater than or equal to a sum of lengths of proper intervals generated in MPA:
$$\sum_{j\in J(\mathrm{PA})} |a_j - a_{j-1}| \geqslant \sum_{j\in J(\mathrm{MPA})} |a_j - a_{j-1}|,$$
where $J(\mathrm{PA})$ and $J(\mathrm{MPA})$ are sets of indices of proper intervals in Petras' algorithm and MPA, respectively; it is obvious that for bad intervals the converse inequality holds;

\item a sum of proper intervals generated in Petras' algorithm is greater than or equal to a sum of proper intervals generated in MPA:
$$\bigcup_{j\in J(\mathrm{MPA})} [a_{j-1},a_j]\subseteq \bigcup_{j\in J(\mathrm{PA})} [a_{j-1}, a_j];$$

\item if MPA stops then Petras' algorithm stops as well, since by (i) we have
$$\sum_{j\not\in J(\mathrm{PA})} |a_j - a_{j-1}| \leqslant \sum_{j\not\in J(\mathrm{MPA})} |a_j - a_{j-1}| \leqslant \frac{\eps}{2M}.$$

\end{enumerate}

Thus we have proved the following lemma.

\begin{lemma}
MPA provides an upper bound for the complexity of Petras' algorithm.
\end{lemma}

\subsection{Leftward flow}\label{subsec:leftward-flow}

The leftward flow is the `opposite' of the rightward flow, i.e.\ we take the longest proper interval in $\D^p_\gamma$, but we move from right to left. Thus we can describe the leftward flow analogously as before:

\begin{enumerate}
\item We start from $y_0 = \beta$ (somewhere in the middle of $(s, s')$, for $s, s'\in S$; the exact value of $\beta$ is irrelevant).

\item Each new interval $[y_{i+1}, y_i]$ is created so that
$\ro(y_{i+1}, y_{i}) \subset \D^p_{\gamma}$ and upper left corner of $\ro(y_{i+1}, y_{i})$  is on the boundary of  $\D^p_{\gamma}$. Observe that all longer intervals $[y_{i+1}-\delta, y_i]$, for $\delta>0$, lack this property.
 Since $\ro(y_{i+1}, y_i) \subset \D_{\gamma}^p$, the interval $[y_{i+1}, y_i]$ must have been accepted during the execution of the MPA and Petras' algorithms.

\item We count the number of steps needed to exit the interval $[\alpha_L, \beta]$ to the left, i.e.\ to have $y_m \leqslant \alpha_L < y_{m-1}$.
\end{enumerate}

This process
creates a sequence of intervals
\begin{equation}\label{left-seq}
y_m \leqslant \alpha_L < y_{m-1} < \ldots < y_0 = \beta.
\end{equation}
Each successive point $y_{n+1}$ is chosen as
\begin{equation}\label{eqn:left-flow}
y_{n+1} = y_n - d_R(y_n),
\end{equation}
where function $d_R(x)$ (depends on the right end of an interval) is equal to $d(x)$, which solves equation (\ref{eqn:next-point1}) with
\begin{equation}\label{eq:gh-leftflow}
h=\frac{B}{2\gamma}, \qquad g=g_R=\frac{A+1}{2}.
\end{equation}

Estimates for $d_R(y)$ are given by Theorem~\ref{thm:d(x)estm-p>1} (for $p>1$), Theorem~\ref{thm:d(x)estm-p<1} (for $p<1$) and (\ref{eq:d-p=1}) (for $p=1$).

\subsubsection{Estimations of $\alpha_L$ for MPA algorithm}

In this subsection we will write $x$ for both the coordinate $x$ and the distance $\dist(x, S)$, and the precise meaning should be clear from the context.
Our goal is to find a lower bound for $\alpha_L$, the distance to the singular point from the right of the set of proper intervals obtained in the MPA algorithm.

Consider two cases:
\begin{enumerate}
\item For $p > 1$ the shape of $\D_{\gamma, S}^p$ is as in Figure~\ref{fig2}.
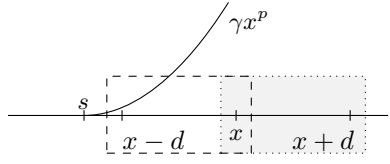
\begin{figure}[htb]
\centering
\begin{tikzpicture}
\draw (1,0) parabola (2.9,1.5)node[anchor=north] {\small \hspace*{4ex}$\gamma x^p$};
\draw (1 cm,2pt) -- (1 cm,-2pt) node[anchor=south] {\small $s_{}$};

\draw[dotted, fill = lightgray!20!white, draw =black] (2.8,-0.5) rectangle (4.7,0.52);
\draw[dashed] (1.3,-0.5) rectangle (3.2,0.52);
\draw (0,0) -- (5,0);

\draw (3.0,2pt) -- (3.0,-2pt) node[anchor=north] {\small $x$};
\draw (4.5,2pt) -- (4.5,-2pt) node[anchor=north east] {\mbox{}\hspace*{1ex} $x+d$\!\!\!};
\draw (1.5,2pt) -- (1.5,-2pt) node[anchor=north west] {\!\!\!\! $x-d$};

\end{tikzpicture}
\caption{An interval $[x, x+d]$ is proper, but $[x-d, x]$ is bad.}\label{fig2}
\end{figure}

 From Theorem~\ref{thm:d(x)estm-p>1} we know that:
\begin{itemize}
\item[] interval $[x, x+d]$  is proper, iff   $d\leqslant d_L(x)$, where $d_L(x) \leqslant c_{1,L}\cdot x^p$,
\item[] interval $[x-d, x]$ is proper, iff  $d\leqslant d_R(x)$, where $c_{2,R}\cdot x^p\leqslant d_R(x)$.
\end{itemize}

Consider $x$ (as in Figure~\ref{fig2}).
Since $[x, x+ d]$ is proper (i.e.\ $d \leqslant d_L(x)$) and $[x-d, x]$ is bad (i.e.\ $d > d_R(x)$), we have:
\begin{equation}\label{eq:d-x-close-to-S}
c_{2,R} \cdot x^p < d \leqslant c_{1,L}\cdot x^p.
\end{equation}
Hence, at any stage of the algorithm, all points $x$ which are the closest to the points from $S$ satisfy the estimate (recall that all bad intervals have the same length):
\begin{equation*}
\left(\frac{d}{c_{1,L}}\right)^{\frac{1}{p}} \leqslant x < \left(\frac{d}{c_{2,R}}\right)^{\frac{1}{p}}.
\end{equation*}
Thus, if
$x_1$ and $x_2$
are the closest to some singular point, then
\begin{equation*}
0 < \left(\frac{c_{2,R}}{c_{1,L}}\right)^{\frac{1}{p}} < \frac{x_1}{x_2} < \frac{\left(\frac{d}{c_{2,R}}\right)^{\frac{1}{p}}}{\left(\frac{d}{c_{1,L}}\right)^{\frac{1}{p}}} = \left(\frac{c_{1,L}}{c_{2,R}}\right)^{\frac{1}{p}} < +\infty.
\end{equation*}

\item For $p \leqslant 1$ from Theorem~\ref{thm:d(x)estm-p<1} we know that for sufficiently small $x$
\begin{itemize}
\item[] interval $[x, x+d]$ is proper, iff $d\leqslant d_L(x)$, where $ d_L(x) \leqslant \frac{x}{g_L}$,
\item[]  interval $[x-d, x]$ is proper, iff $d\leqslant d_R(x)$, where there exists $g_R^+$ such that
$$\frac{x}{g_R^+} < x\left(\frac{1}{g_R} - c_{1,R}x^{\frac{1}{p} - 1}\right) \leqslant d_R(x).$$
\end{itemize}
Since $[x, x+d]$ is proper and $[x-d, x]$ is bad, we have:
\begin{equation}\label{eq:d-x-p<1}
\frac{x}{g_R^+} < d \leqslant \frac{x}{g_L}.
\end{equation}
Hence, at any stage of the algorithm, all points $x$ which are the closest to the points from $S$ satisfy the estimate (recall that all bad intervals have the same length):
\begin{equation*}
d\cdot g_L \leqslant x < d\cdot g_R^+.
\end{equation*}
Thus if we have points $x_1$ and $x_2$ which are the closest to some singular point, then
\begin{equation*}
0 < \frac{A - 1}{A +1} \approx \frac{g_L}{g_R^+} \leqslant \frac{x_1}{x_2} \leqslant \frac{d\cdot g_R^+}{d\cdot g_L} = \frac{g_R^+}{g_L} \approx \frac{A + 1}{A - 1} < +\infty.
\end{equation*}

\end{enumerate}

\begin{figure}[htb]
\centering
\begin{picture}(200,50)
\put(0,25){\line(1,0){280}}
\put(10,20){\line(0,1){10}}
\put(2,10){$-1$}

\put(50,23){\line(0,1){25}}
\put(46,17){$x_1$}
\put(65,23){\line(0,1){25}}
\put(63,17){$x_2$}
\put(50,48){\line(1,0){15}}
\put(60,20){\line(0,1){10}}
\put(55,35){$s_1$}

\put(90,23){\line(0,1){25}}
\put(86,17){$x_3$}
\put(120,23){\line(0,1){25}}
\put(115,17){$x_4$}
\put(90,48){\line(1,0){30}}
\put(110,20){\line(0,1){10}}
\put(105,35){$s_2$}

\put(150,20){\line(0,1){10}}
\put(160,15){$\ldots$}
\put(180,20){\line(0,1){10}}

\put(210,23){\line(0,1){25}}
\put(198,17){$x_{2m - 1}$}
\put(240,23){\line(0,1){25}}
\put(210,48){\line(1,0){30}}
\put(240,20){\line(0,1){10}}
\put(237,14){$s_m = 1$}

\end{picture}
\caption{Singularities $s_1, \ldots, s_m$ and blocks of bad intervals between $x_{2i-1}$ and $x_{2i}$ for $i=1,\ldots, m$.}\label{fig:bad-int}
\end{figure}
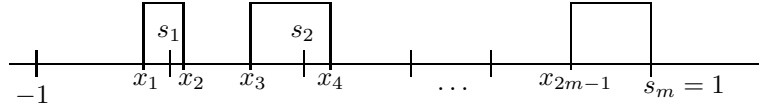

Now we know that for any $p>0$, at any stage of the MPA, for any points $x_i$, $x_j$ which are the closest to a singular point from $S$:
\begin{equation}\label{eqn:dist}
0 < \frac{x_i}{x_j} \leqslant T < +\infty,
\end{equation}
where $T$ is a constant independent of $\eps$. Because~(\ref{eqn:dist}) holds for the minimal and maximal distance denoted by $x_{\min}$, $x_{\max}$ we have
\begin{equation}\label{eqn:dist-min-max}
x_{\max}\leqslant T\cdot x_{\min}.
\end{equation}
Bad intervals are in the neighbourhood of any singular point from $S$, so there exists $I\leqslant 2m$ such that for any $1\leqslant i\leqslant I$ a point $x_i$ is an end of the connected block of proper intervals (see Figure~\ref{fig:bad-int}).

While the algorithm is running
\begin{equation*}
\frac{\eps}{2M} < \sum_{1\leqslant i\leqslant I} x_i \leqslant \sum_{1\leqslant i\leqslant I} x_{\max}
\end{equation*}
and from~(\ref{eqn:dist-min-max}) we have
\begin{equation*}
\frac{\eps}{2M} < \sum_{1\leqslant i\leqslant I} T\cdot x_{\min} \leqslant 2m\cdot T\cdot x_{\min}
\end{equation*}
obtaining
\begin{equation}\label{eq:xmin-estm}
x_{\min} \geqslant \frac{\eps}{4mMT}.
\end{equation}
Let us stress that (\ref{eq:xmin-estm}) holds as long as the stopping condition in MPA algorithm is not satisfied. We need to estimate how far we can go
in the last stage (i.e.\ $\alpha_L$).

Let $x_1$ be a point from the proper interval, which is the closest to some $s \in S$ from right or left. For simplicity we will assume that $s <x_1$, the other case is analogous. Let $d$ be the length
of bad intervals, which will be now divided by $2$. Let $x_2$ (with $s <x_2 < x_1$) be such that $[x_2,x_1]$ is covered by proper intervals of length $d/2$.  Therefore interval $[x_2,x_2 + d/2]$ is proper.

\begin{enumerate}
\item Case $p>1$.
By~(\ref{eq:d-x-close-to-S}),~(\ref{eq:xmin-estm}) and Theorem~\ref{thm:d(x)estm-p>1} we have the following estimations
$$
\left.
\begin{array}{lclcl}
 & & c_{2,R}\cdot x_1^p & < & d\\
\displaystyle \frac{d}{2} & < & d_L(x_2) & \leqslant & c_{1,L}\cdot x_2^p
\end{array}
\right\} \Longrightarrow\quad   x_2\geqslant \left(\frac{1}{2}\cdot \frac{c_{2,R}}{c_{1,L}}\right)^{1/p} x_1 \geqslant
\left(\frac{1}{2}\cdot\frac{c_{2,R}}{c_{1,L}}\right)^{1/p} \frac{\eps}{4mMT}.
$$

\item Case $p\leqslant 1$. By (\ref{eq:d-x-p<1}),~(\ref{eq:xmin-estm}) and Theorem~\ref{thm:d(x)estm-p<1} we obtain
$$
\left.
\begin{array}{lclcl}
& & \displaystyle\frac{x_1}{g_R^+} & < & d\\
\displaystyle\frac{d}{2} & < & d_L(x_2)  & \leqslant & \displaystyle\frac{x_2}{g_L}
\end{array}
\right\} \Longrightarrow\quad x_2\geqslant \frac{1}{2}\cdot\frac{g_L}{g_R^+}\cdot x_1 \geqslant  \frac{1}{2}\cdot\frac{g_L}{g_R^+}\cdot \frac{\eps}{4mMT}.
$$

\end{enumerate}
Observe that in both above cases we have obtained
\begin{equation}
  \alpha_L \geqslant \frac{\eps}{M T_0}, \label{eq:estm-alpha0}
\end{equation}
for some constant $T_0$.

\subsection{Estimation of the number of proper intervals from above}

\subsubsection{Case $p>1$}
\begin{lemma}\label{lm:upper-bound}
Assume that $p > 1$ and $f$ satisfies \PPCS{p, \gamma, \{0\}}. Then
$$\mathcal{Z}(\D, \eps) = {\cal O}\left(\left(\frac{1}{\eps} \right)^{p-1} \right).$$
\end{lemma}

\begin{proof}
Let us assume that $[\alpha_i, \alpha_{i+1}]$ is obtained during MPA by bisecting an interval $[u,v]$. Since $[u,v]$ was bisected, it was bad, thus
\begin{equation}
v - u > d_R(v).  \label{eq:v-u-dl}
\end{equation}
 From  Theorem~\ref{thm:d(x)estm-p>1} we know that $d_R(x)\geqslant c_2\cdot x^p$, where $c_2 >0$. Let us denote $$w(x) = c_2\cdot x^p.$$

\begin{figure}[ht]
\centering
\begin{minipage}[t]{0.4\linewidth}
\begin{picture}(150,50)
\put(0,40){\line(1,0){150}}

\put(10,38){\line(0,1){4}}
\put(45,30){\line(0,1){12}}
\put(85,30){\line(0,1){12}}
\put(125,30){\line(0,1){12}}

\put(8,28){0}
\put(42,24){$\alpha_i$}
\put(82,23){$\alpha_{i+1}$}
\put(122,23){$z$}
\put(122,12){$\parallel$}
\put(50,1){$\alpha_{i+1} + (\alpha_{i+1} - \alpha_i)$}
\put(60,45){L}
\put(100,45){R}
\end{picture}
\caption{The case (L)}\label{fig:L}
\end{minipage}
\quad
\begin{minipage}[t]{0.4\linewidth}
\begin{picture}(150,50)
\put(0,40){\line(1,0){150}}

\put(10,38){\line(0,1){4}}
\put(45,30){\line(0,1){12}}
\put(85,30){\line(0,1){12}}
\put(125,30){\line(0,1){12}}

\put(8,28){0}
\put(42,24){$z$}
\put(42,12){$\parallel$}
\put(42,1){$\alpha_i - (\alpha_{i+1} - \alpha_i)$}
\put(82,23){$\alpha_i$}
\put(122,23){$\alpha_{i+1}$}
\put(60,45){L}
\put(100,45){R}
\end{picture}
\caption{The case (R)}\label{fig:R}
\end{minipage}
\end{figure}

Two cases are possible as a result of the bisection:
\begin{enumerate}
\item[(L)] $[\alpha_i, \alpha_{i+1}]$ is the left part  in $[\alpha_i,z]$ (see Figure~\ref{fig:L}), where
$z = \alpha_{i+1}+(\alpha_{i+1} - \alpha_i)$, and from (\ref{eq:v-u-dl}) we have
\begin{eqnarray}
z - \alpha_i = 2(\alpha_{i+1} - \alpha_i) & > & d_R(\alpha_{i+1} + (\alpha_{i+1} - \alpha_i))\nonumber\\
  & \geqslant & w(\alpha_{i+1} + (\alpha_{i+1} - \alpha_i)) > w(\alpha_{i+1}). \nonumber
\end{eqnarray}

\item[(R)] $[\alpha_i, \alpha_{i+1}]$ is the right part in $[z,\alpha_{i+1}]$ (see Figure~\ref{fig:R})  so for both $z > 0$ and $ z\leqslant 0$ it holds
$$
\begin{array}{lclp{3cm}}
 2(\alpha_{i+1} - \alpha_i)=\alpha_{i+1} - z & > & d_R(\alpha_{i+1}), & by~(\ref{eq:v-u-dl})\\
  & \geqslant & c_2\cdot \alpha_{i+1}^p\\
  & = & w(\alpha_{i+1}).
\end{array}
$$
\end{enumerate}

Therefore in both cases the following estimate holds:
\begin{equation}\label{eqn:estimation-below}
\alpha_{i+1} - \alpha_i > \frac{1}{2}\cdot w(\alpha_{i+1}).
\end{equation}

In the leftward process we are moving from right to left thus we have to look at these $\alpha$'s in the opposite direction:
$$
\begin{array}{rlllll}
\displaystyle\frac{\eps}{M T_0} = \alpha_0 & \ldots & \alpha_i & \alpha_{i+1} & \ldots & \alpha_n = \beta\\
y_n & \ldots & y_{i+1} & y_i & \ldots & y_0.
\end{array}
$$
Consequently substituting $\alpha_{n - i}$ by $y_i$ ($0\leqslant i \leqslant n$) in~(\ref{eqn:estimation-below}) we obtain exactly the recurrence of the leftward flow:
\begin{eqnarray}\label{eqn:leftward-final}
y_{i+1} & \leqslant & y_i - \frac{1}{2}\cdot \underbrace{c_2\cdot y_i^p}_{w(y_i)}.
\end{eqnarray}
It is easy to see that
\begin{equation*}
 y_i \leqslant y(i),
\end{equation*}
where $y(t)$ is a solution of
\begin{equation}\label{eqn:diff}
y'= - a\cdot y_i^p,
\end{equation}
with $\displaystyle a = \frac{1}{2}\cdot c_2$ and $y(0)=y_0$.
Knowing that the solution of~(\ref{eqn:diff}) is of the form
\begin{equation*}
  y(t)=\frac{y(0)}{\left(1 + a\cdot y(0)^{p-1}(p-1)t \right)^{\frac{1}{p-1}}}
\end{equation*}
and setting
$y(0) = \beta$ and $\displaystyle y(t) = \frac{\eps}{M T_0}$
(compare (\ref{eq:estm-alpha0})) we calculate
\begin{eqnarray*}
t & = & \frac{\left(\frac{M T_0 \beta}{\eps}\right)^{p-1}-1}{a \beta^{p-1}(p-1)}\\
  & = & \frac{(M T_0)^{p-1}}{\eps^{p-1}}\cdot\frac{1}{a(p-1)} - \frac{1}{a(p-1)\beta^{p-1}}.
\end{eqnarray*}
Observe that $t$ is un upper bound for $n$,
thus (see Definition~\ref{def:asym-rate-groth})
\begin{equation*}
\mathcal{Z}(\D, \eps) =  {\cal O}\left(\left(\frac{1}{\eps} \right)^{p-1} \right).
\end{equation*}
\qed
\end{proof}

The above theorem holds also for $|S| = m\geqslant 1$, but in such a case $\beta$ is in the middle of $[s_{i-1}, s_i]$ for each $1 < i \leqslant m$.

\subsubsection{Case $p \leqslant 1$}

\begin{lemma}\label{lm:upper-estimation-p<1}
Assume that $p \leqslant 1$ and $f$ satisfies \PPCS{p, \gamma, \{0\}}. Then there exists $\eta > \eps$ such that in the segment $[0, \eta]$
$$\mathcal{Z}(\D, \eps) = \bigo\left(\ln \frac{1}{\eps} \right).$$
\end{lemma}

\begin{proof}
Following~(\ref{eq:gh-leftflow}) we set
\begin{equation*}
  g_R=\frac{A+1}{2} >1. 
\end{equation*}

Let $c_1$ be as in Theorem~\ref{thm:d(x)estm-p<1}.
Let us take
\begin{equation*}
\eta = \frac{1}{(2g_R c_1)^{\frac{p}{p-1}}}
\end{equation*}
and notice that for all $x\in [0,\eta]$
\begin{equation*}
\frac{1}{g_R} - c_1\cdot x^{\frac{1}{p}-1} \geqslant \frac{1}{2g_R}
\end{equation*}
thus from Theorem~\ref{thm:d(x)estm-p<1} we have
\begin{equation}
\frac{x}{2g_R} \leqslant d_R(x).  \label{eq:estm-d2}
\end{equation}

Let $k$ be the number of proper intervals in $[0,\eta]$.
Now, as in Lemma~\ref{lm:upper-bound}, we have two cases.
In the case (R) (see Figure~\ref{fig:R}) we have
\begin{eqnarray*}
\alpha_{i+1} - z = 2(\alpha_{i+1} - \alpha_i)  \geqslant  d_R(\alpha_{i+1}) \geqslant  \frac{\alpha_{i+1}}{2g_R}.
\end{eqnarray*}

In the case (L) (see Figure~\ref{fig:L}) using (\ref{eq:estm-d2}) we obtain
\begin{eqnarray*}
z - \alpha_i = 2(\alpha_{i+1} - \alpha_i) & \geqslant & d_R(z) \geqslant \frac{z}{2g_R}\\
  & \geqslant & \frac{\alpha_{i+1}+(\alpha_{i+1} - \alpha_i)}{2g_R}  \geqslant  \frac{\alpha_{i+1}}{2g_R}.
\end{eqnarray*}
In both cases we get
\begin{equation*}
  \left(1 - \frac{4}{g_R}\right) \alpha_{i+1} > \alpha_i.
\end{equation*}
Since $g_R>1$, we obtain
\begin{equation*}
\alpha_{i+1} \geqslant \alpha_i\left(1 + \frac{1}{4g_R - 1} \right).
\end{equation*}

Therefore from  (\ref{eq:estm-alpha0})
\begin{equation}\label{eqn:t-upper-p<1}
 \eta \geqslant \alpha_{i} \geqslant \alpha_0\left(1 + \frac{1}{4g_L - 1} \right)^{i} \geqslant \frac{\eps}{M T_0}\left(1 + \frac{1}{4g_L - 1} \right)^{i}.
\end{equation}

Hence we obtain the following upper bound for $k$
\begin{equation*}
  k \leqslant  \frac{\ln \left(\frac{\eta M T_0}{\eps}\right)}{\ln \left( 1+ \frac{1}{4g_L - 1}\right)} +1.
\end{equation*}
\qed
\end{proof}

\section{Complexity of Petras' algorithm}
\label{sec:complexity}

Finally we can state a theorem estimating the complexity of Petras' algorithm for functions satisfying PPC and NPC conditions. As we mentioned in the introduction by the complexity we understand the number of evaluations of an integrand at the nodes produced by the algorithm. Therefore  this is not a comprehensive evaluation, because we  neglect  complexity of checking analyticity, calculating bounds and the precision of arithmetic operations.

\begin{theorem}
\label{thm:main-complexity}
Assume that $S = \{s_1, \ldots, s_m\}$, $m \geqslant 1$.
For functions satisfying \PPCS{p, \gamma, S} and \NPC{p, \gamma', s_0, \beta} for some $s_0 \in S$, the complexity of Petras' algorithm is
\begin{eqnarray}
\Theta(|\ln \eps|/\eps^{p-1}), & & \mbox{for } p > 1,\label{eqn:p>1}\\
\Theta(\ln^2 \eps), & & \mbox{for } p \leqslant 1.\label{eqn:p<1}
\end{eqnarray}

\end{theorem}

\begin{proof}
First notice that by Lemmas~\ref{lm:lower-bound} and~\ref{lm:n-estimation-p<1} for  $s_0\in S$
the number of proper intervals created during the execution of the algorithm is
\begin{eqnarray}
\Omega(1/\eps^{p-1}), & & \mbox{for } p > 1,\\
\Omega(|\ln \eps|), & & \mbox{for } p \leqslant 1.
\end{eqnarray}

On the other hand for some small neighbourhood of each $s \in S$ in Lemmas~\ref{lm:upper-bound} and~\ref{lm:upper-estimation-p<1} we obtained an upper bound for
 the number of proper intervals created during the execution of the algorithm:
\begin{eqnarray}
\bigo(1/\eps^{p-1}), & & \mbox{for } p > 1,\label{eqn:upp-p>1}\\
\bigo(|\ln \eps|), & & \mbox{for } p \leqslant 1\label{eqn:upp-p<1}.
\end{eqnarray}

Next use Lemma~\ref{lm:analytic-fun} to estimate the number of proper intervals outside the neighbourhood of $S$. This estimate is a constant (dependent on $f$).

Hence we obtain the assertion of the theorem.
\qed
\end{proof}

\section{Functions giving rise to the complexity $\bigo\left(\frac{|\ln\eps|}{\eps^{p-1}}\right)$}\label{sec:examples}
In Section~\ref{subsec:petras-cond} we have shown that function $f(z) = \sin(1/z)$ satisfies \PPC{2, \gamma, \{0\}} and \NPC{2, \gamma, 0, 1}. The complexity of Petras' algorithm for such functions is $\displaystyle \bigo\left(\frac{|\ln\eps|}{\eps}\right)$. In this section we give examples of functions for which the complexity is $\displaystyle \bigo\left(\frac{|\ln\eps|}{\eps^{p-1}}\right)$, for any $p > 1$.

\begin{theorem}
\label{thm:petras-ppc-p}
Let $S = \{0\}$ and let $f(z) = \sin (1/z^{p-1})$, where $p \in \nat$, $p > 1$.
Then for any $c > 1$
there exists $\gamma>0$ such that
\begin{equation*}
  \sup_{z \in D^p_{\gamma} } |f(z)| \leqslant c \cdot \sup_{x \in [-1,1]}|f(x)|.
\end{equation*}

Therefore  for any  $\alpha, \beta\in [-1,1]$
\begin{equation*}
\ro(\alpha, \beta)\subset \D^p_{\gamma} \quad\Longrightarrow\quad f \mbox{ is analytic on } \ro(\alpha, \beta)\ \wedge\ \sup_{z\in\ro(\alpha, \beta)}\left|f(z)\right|\leqslant c \cdot \sup_{x \in [-1,1]} |f(x)|.\\
\end{equation*}

\end{theorem}
\begin{proof}
Let $\gamma \leqslant 1$.
If $\ro(\alpha, \beta)\subset \D^p_{\gamma}$ then it is obvious that $f$ is analytic on $\ro(\alpha, \beta)$. Thus it is enough to check the second part of the thesis.

Let us set $n=p-1$.  Assume $z=x+iy \in D^p_{\gamma} \cap \{ z: |\re{z}|  \leqslant 1\}$. Let  $r^2 = x^2 + y^2$.
As in the proof of Theorem~\ref{thm:petras-insufficient}, we have
\begin{eqnarray*}
  |\sin(1/z^n) | \leqslant \frac{1}{2} \left(\exp\left(|\im{1/z^n}|\right) + \exp\left(-|\im{1/z^n}|\right)\right).
\end{eqnarray*}

By~(\ref{eqn:1-z}) we have
\begin{eqnarray*}
  \im{1/z^n} = \frac{1}{r^{2n}} \left( \binom{n}{1}x^{n-1}y  - \binom{n}{3}x^{n-3}y^3 + \dots \right)
\end{eqnarray*}
and it is easy to see that for $z \in D^p_{\gamma} \cap \{ z: |\re{z}| \leqslant 1  \}$ and for some constant $\Phi=\Phi(n)$,
\begin{eqnarray*}
 \left|\im{1/z^n}\right| \leqslant  \frac{\Phi}{r^{2n}}  |x^{n-1} y| \leqslant \Phi \frac{\gamma |x|^{n-1} |x|^p }{|x|^{2n}}= \Phi \gamma.
\end{eqnarray*}

For $z \in D^p_\gamma \cap \{z: |\re{z}| \leqslant 1\}$ we have (because $x \mapsto x + \frac{1}{x}$ is increasing for $x >1$)
\begin{equation*}
   |\sin(1/z^n) | \leqslant \frac{1}{2} \left(\exp(\Phi\gamma) + \exp(-\Phi\gamma) \right) \to 1, \quad \mbox{as } \gamma\to 0.
\end{equation*}
Thus there exists $c>1$ such that taking $\gamma$ sufficiently small we get
\begin{equation*}
\sup_{z \in D^p_{\gamma}\cap \{z:|\re{z}| \leqslant 1\}} \left| \sin \left(1/z^{p-1} \right)\right|
< c \cdot \sup_{x \in [-1,1]} |\sin(1/x^{p-1})|.
\end{equation*}

\begin{figure}[htb]
\centering

\begin{tikzpicture}
\draw[draw=black ] (1,0) parabola (3.5,1.15);
\filldraw[fill=gray!20!white, draw=gray!20!white ] (1,0) parabola (3.5,1.15) -- (3.5,0);
\filldraw[fill=gray!20!white, draw=gray!20!white ] (1,0) parabola (3.5,-1.15) -- (3.5,0);
\draw (1 cm,1pt) -- (1 cm,-1pt) node[anchor=north] {\small $0$};
\draw (3.5 cm,1pt) -- (3.5 cm,-1pt) node[anchor=north] {\small $1$};
\draw[draw=black ] (1,0) parabola (3.5,1.15);
\draw[draw=black ] (1,0) parabola (3.5,-1.15);

\draw[dashed] (2.1,-0.22) rectangle (3.8,0.22);

\draw (3.5,0.22) -- (3.8,0.22);
\draw (3.8,0) -- (3.8,0.22);

\draw[->] (3.8, 0.7) -- (3.7, 0.23);
\draw[] (3.8, 0.65) -- (3.8, 0.65) node[anchor=south] {\tiny $dx(\ro)$};
\draw[->] (4.3, 0.15) -- (3.8, 0.1);
\draw[] (4.2, 0.15) -- (4.2, 0.15) node[anchor=west] {\tiny $dy(\ro)$};

\draw[color = black] (0,0) -- (5,0);
\end{tikzpicture}
\caption{Rectangle $\ro(\alpha, \beta)$ does not have to be entirely contained in $D^p_{\gamma}\cap \{z: |\re{z}| \leqslant 1\}$, but for small $\gamma$ the projecting part ($\ro(\alpha,\beta) \backslash (D^p_{\gamma}\cap \{z: |\re{z}| \leqslant 1\}) = 2(dx(\ro)\times dy(\ro))$) is small.}\label{fig:proj-part}
\end{figure}

Note that (see Figure~\ref{fig:proj-part}) the rectangle $\ro(\alpha,\beta)$ does not have to be entirely contained in the region $D^p_{\gamma}\cap \{z: |\re{z}| \leqslant 1\}$, but the width and height of the projecting area linearly depend on $\gamma$ ($\displaystyle dx(\ro) < \gamma \frac{A-1}{B}$ and $dy(\ro) < \gamma$ --- see~Theorem \ref{thm:G-d-exists}) and the function is continuous in the wide neighbourhood of 1, therefore further reducing $\gamma$ we obtain
\begin{equation*}
\sup_{z\in\ro(\alpha,\beta)} \left|\sin \left(1/z^{p-1} \right)\right| \leqslant c \cdot \sup_{x \in [-1,1]} |\sin(1/x^{p-1})|.
\end{equation*}
\qed
\end{proof}

Theorem~\ref{thm:petras-ppc-p} says that for $z\mapsto \sin(1/z^{p-1})$ there exists a region $\D^p_{\gamma} \cap \{z: |\re{z}| \leqslant 1\}$ where the function is analytic and appropriately bounded. The next theorem says that
there exists a region of the same shape as before, but such that (on the boundary of this region
close to the singular point) the values of the function are arbitrarily large.

\begin{theorem}
\label{thm:petras-condition-3}Let $S = \{0\}$ and let $f(z) = \sin(1/z^{p-1})$, where $p\in \nat$, $p > 1$.
Then for any $c > 1$ there exist $\gamma > 0$ and $\eta >0$ such that for any $\alpha\in [-\eta,\eta]$ and $\beta\in (\alpha, 1]$
\begin{equation*}
\ro(\alpha, \beta)\not\subset \D^{p}_{\gamma} \quad\Longrightarrow\quad f \mbox{ is not analytic on } \ro(\alpha, \beta)\ \vee \sup_{z\in\ro(\alpha, \beta)}\left|f(z)\right| > c \sup_{x \in [-1,1]} |f(x)|.
\end{equation*}
\end{theorem}
\begin{proof}
Let us fix $\gamma >0$ and consider points $z=x+iy \in \partial D^p_{\gamma} \cap \{z:\ |\re{z}| \leqslant 1\}$.  We will restrict our attention to $x>0$ and $y=\gamma x^{p} > 0$.
By~(\ref{eqn:1-z}), for sufficiently small $x$ and for some constants $\Phi, \Phi'>0$ there is 
\begin{eqnarray*}
  \im{\frac{1}{z^{p-1}}} & = & \frac{1}{r^{2(p-1)}} \left( \binom{p-1}{1}x^{p-2}y - \binom{p-1}{3}x^{p-4}y^3 + \dots \right)\\
  & \geqslant &
    \frac{\gamma\Phi x^{2(p-1)}}{r^{2(p-1)}} \\
  & = &      \frac{\gamma\Phi x^{2(p-1)}}{x^{2(p-1)}(1 + \gamma^2 x^{2(p-1)})^{p-1}}\\
  & = & \frac{\gamma\Phi}{(1 + \gamma^2 x^{2(p-1)})^{p-1}}\\
  & \geqslant & \gamma\Phi'.
\end{eqnarray*}

Observe that (compare the proof of Theorem~\ref{thm:petras-condition-2}), since $x \mapsto x - \frac{1}{x}$, we have
\begin{eqnarray*}
\inf_{z\in\W_{\gamma}}  |\sin (1/z^l)| & \geqslant & \frac{1}{2} \left(\exp( |\im{1/z^l}|) - \exp( -|\im{1/z^l}|) \right)\\
& \geqslant & \frac{1}{2} \left(\exp(\gamma\Phi') - \exp(-\gamma\Phi') \right),
\end{eqnarray*}
where $\W_\gamma = \partial D^p_{\gamma} \cap \{z: |\re{z}| \leqslant \eta\}\backslash \{0\}$.

Now, for any $c > 1$ it is enough to take $\gamma > 0$ big enough to have $\inf_{z\in\W_{\gamma}}  |\sin (1/z^l)| > c$.
\qed
\end{proof}

From the above theorems it follows  that functions $f(z)= \sin (1/z^{p-1})$ for $p \in \nat$, $p > 1$ satisfy the assumptions of Theorem~\ref{thm:main-complexity}, thus their complexity scales as $|\ln \eps|/\eps^{p-1}$.

\section{Appendix: asymptotic rate of growth}\label{apx:asym-rate-groth}

Let us adopt the notation used in the classical (discrete) complexity theory to describe asymptotic rate of growth.

\begin{definition}\label{def:asym-rate-groth}
Let $f, g:\reals\to \reals_+$ .
\begin{enumerate}
\item We say that $f$ {\em is at least of order} $g$, if there exist $\eps_0 >0$ and $c>0$, such that:

$$
\forall\ 0 < x \leqslant \eps_0 :\ f(x) \geqslant c \cdot g(x).
$$
Notation: $f(x) \in \Omega(g(x))$.

\item We say that $f$ {\em is at most of order} $g$, if there exist $\eps_0 >0$ and $c>0$, such that:

$$
\forall\ 0 < x \leqslant \eps_0 :\ f(x) \leqslant c \cdot g(x).
$$
Notation: $f(x) \in {\cal O}(g(x))$.

\item We say that $f$ {\em is exactly of order} $g$, if there exist $\eps_0 >0$, $c_1>0$ and $c_2$, such that:

$$
\forall\ 0 < x \leqslant \eps_0 :\ c_1 \cdot g(x) \leqslant f(x) \leqslant c_2 \cdot g(x).
$$
Notation: $f(x) \in \Theta(g(x))$.

\end{enumerate}
\end{definition}

\bibliographystyle{plain}
\bibliography{ref}

\begin{thebibliography}{1}

\bibitem{Br}
H.~Brass.
\newblock {\em Quadraturverfahren}.
\newblock Vandenhoeck and Ruprecht, G{\"o}ttingen, 1997.

\bibitem{GaertnerHotz12}
T.~Gaertner and G.~Hotz.
\newblock Representation theorems for analytic machines and computability of
  analytic functions.
\newblock {\em Theory of Computing Systems}, 51:65--84, 2012.

\bibitem{Hoeven05}
J.~{\VAN{Hoeven}{Van}{van der}}~Hoeven.
\newblock Effective analytic functions.
\newblock {\em Journal of Symbolic Computation}, 39:433--449, 2005.

\bibitem{P98}
K.~Petras.
\newblock Gaussian {Versus} {Optimal} {Integration} of {Analytic} {Functions}.
\newblock {\em Constructive Approximation}, 14:231--245, 1998.

\bibitem{P98-2}
K.~Petras.
\newblock On the complexity of self-validating numerical integration and
  approximation of functions with singularities.
\newblock {\em Journal of Complexity}, 14:302--318, 1998.

\bibitem{P02}
K.~Petras.
\newblock Self-validating integration and approximation of piecewise analytic
  functions.
\newblock {\em Journal of Computational and Applied Mathematics}, 145:345--359,
  2002.

\end{thebibliography}

\end{document}